\documentclass[11pt]{article}

\usepackage[margin=1in]{geometry}
\usepackage{amsmath,amssymb,amsthm,mathtools}
\usepackage{bm}
\usepackage{graphicx}
\usepackage{xcolor}
\usepackage{enumitem}
\usepackage[round]{natbib}
\usepackage[colorlinks=true,linkcolor=blue,citecolor=blue,urlcolor=blue]{hyperref}

\theoremstyle{plain}
\newtheorem{theorem}{Theorem}
\newtheorem{lemma}{Lemma}
\newtheorem{proposition}{Proposition}
\newtheorem{corollary}{Corollary}
\theoremstyle{definition}
\newtheorem{definition}{Definition}
\newtheorem{example}{Example}
\newtheorem{remark}{Remark}

\def\reals{\mathbb R}

\def\tTheta{\tilde{\Theta}}
\def\tbX{\tilde{\bm X}}

\title{Precise Asymptotics for Dense High-Dimensional Huber Regression under Heteroscedasticity}

\author{Jing Zhou\thanks{Department of Mathematics, University of Manchester. Email: \texttt{jing.zhou@manchester.ac.uk}.}
\and Ramji Venkataramanan\thanks{Department of Engineering, University of Cambridge. Email: \texttt{rv285@cam.ac.uk}.}
\and Hui Zou\thanks{Department of Applied Mathematics and Statistics, Johns Hopkins University. Email: \texttt{hzou20@jh.edu}.}}

\date{\today}

\begin{document}

\maketitle

\begin{abstract}
We study Huber regression in a heteroscedastic linear model in the proportional regime $n/p \to \delta \in (1, \infty)$. Our main result is an exact characterization of the asymptotic empirical distribution of the Huber estimator. This characterization is in terms of three scalar parameters which determine the first-order centre  and the asymptotic variance of the Huber estimator.  The result shows that with heteroscedastic noise, the first-order centre  may contain an additional component along the variance index coefficient. We derive the displacement coefficients that determine this centre, as well as the exact asymptotic variance. The analysis uses the matrix approximate message passing (AMP) framework to jointly track the regression and variance index coefficients, extending the existing framework for homoscedastic M-estimators. We illustrate the implications of this characterization by constructing a test for centre  displacement and examining its empirical size and power through simulations.
\end{abstract}


\section{Introduction} \label{sec:intro}

Huber regression \citep{Huber1964,huber1973robust,HuberRonchetti2009} is a standard robust alternative to least squares, often used when the data have outliers or the noise is heavy-tailed. Consider a scalar random variable $Y \in \mathbb R$. Its $a$-Huber mean, that is the Huber mean at threshold $a \in (0, \infty] $, is defined as
\begin{equation}\label{eq: population target huber mean}
    m_a (Y) = \arg\min_{z \in \mathbb R} E [\ell_a (z, Y)],
\end{equation}
where
\begin{eqnarray}\label{eq: Huber loss function}
  \ell_a (z, y) =
  \begin{cases}
\frac{1}{2} (z- y) ^2, & |z - y| \leq a \\
a |z - y| - \frac{1}{2} a^2 , & |z - y| > a.
\end{cases}
\end{eqnarray}

We study Huber regression for the following linear model with heteroscedastic noise:
\begin{equation}\label{eq: data generating process}
    Y=X^\top\beta_0+\sigma(X^\top\gamma_0)\bar\varepsilon,
\end{equation}
where $Y, \bar\varepsilon \in \mathbb R$, and $X \in \mathbb R^p$. The noise variable $\bar \varepsilon$ is independent of the covariate vector $X$. The regression parameter $\beta_0 \in \mathbb R^p$ is of primary interest. 
Moreover, the conditional variance of $Y$ depends on $X$ through a variance index parameter vector $\gamma_0 \in \mathbb R^p$ and an (unknown) scale function $\sigma: \mathbb R\to(0,\infty)$.

We consider a sample $\{(Y_i, X_i)\}_{i=1}^n$ from the heteroscedastic model in \eqref{eq: data generating process} such that $Y_i = X_i^\top\beta_0 + \sigma(X_i^\top \gamma_0)\varepsilon_i$.   
The regression parameter $\beta_0$ can be estimated by solving 
\begin{equation}\label{eq: M-estimator minimization}
    \hat\beta\in\arg\min_{\beta\in\mathbb R^p} \sum_{i=1}^n\ell_a(X_i^\top\beta,Y_i),
\end{equation}
where $\ell_a$ is the Huber loss function in \eqref{eq: Huber loss function}.

In classical fixed-dimensional asymptotics, where $p$ is constant and $n\to \infty$, the least-squares estimator remains consistent for $\beta_0$ under a correctly specified linear conditional mean and standard regularity conditions. Heteroscedasticity affects its asymptotic covariance. Inference can therefore be
based on heteroscedasticity-consistent covariance estimation in the Eicker--Huber--White tradition
\citep{Eicker1967,Huber1967,White1980}, while suitable weighting can improve efficiency \citep{carroll2017transfor}.

Relatively little is known about the behaviour of Huber regression under heteroscedasticity in the dense high-dimensional regime. We study the proportional asymptotic setting, where the sample size $n$ and the dimension $p$ both grow, with $n/p\to\delta\in(1,\infty)$ as $p \to \infty$.  Our main result, Theorem~\ref{thm: dist convergence of huber mean}, precisely characterizes the asymptotic empirical distribution of the Huber estimator, showing that the limiting centre  of the Huber estimator admits an explicit decomposition into the regression parameter and variance index parameter directions. This yields an exact expression for the asymptotic variance:
        \begin{equation}
            \label{eq: asymptotic risk convergence init}
            \frac{1}{p} (\hat\beta - \beta_0 - \nu_* \gamma_0)^\top \Omega_p (\hat\beta -  \beta_0 - \nu_* \gamma_0) \stackrel{a.s.}{\to} \tau_*^2,
        \end{equation}
    where $\Omega_p$ is the scaled population covariance matrix, $ \nu_*$ and $\tau_*^2$ are deterministic constants, 
    and $\stackrel{a.s.}{\to}$ denotes almost sure convergence. We validate this result through simulation in Section~\ref{ssec: work out example}. Moreover, we derive an exact effective-score orthogonality criterion for $\nu_* = 0$. In particular, the nuisance coefficient $\nu_*$ vanishes if the derived \emph{effective Huber score} is correctly centred. Under this condition, the limiting centre  contains no displacement along $\gamma_0$, establishing Fisher consistency with respect to the heteroscedastic nuisance direction.

The proof of our main result is based on an iterative matrix-valued  Approximate Message Passing (AMP) algorithm \citep{wang2024universality, tan2024approximate},  tailored to coincide with  the Huber estimator at convergence. The asymptotic characterization of the matrix-valued AMP simultaneously tracks the effect of  $\beta_0$ and $\gamma_0$ on the iterates. This characterization is then translated to the Huber estimator to obtain the main result.  Our proof generalizes the vector-valued AMP framework used by \cite{donoho2016high} for homoscedastic M-estimation. A further technical challenge is that without additional assumptions, the minimizer of \eqref{eq: M-estimator minimization} need not be unique.  
The analysis of homoscedastic M-estimation by \citet{donoho2016high} relies on the uniform strong convexity condition $\inf_x \ell''(x) > 0$. We use a milder assumption, discussed in Section~\ref{sec: Huber mean main section}.

Convex estimators have been analysed in the proportional asymptotic regime for various regression models with homoscedastic noise. We mention the most relevant works here, and refer the reader to \cite{feng2022unifying} for a broader overview. 
For the homoscedastic linear model,  \cite{el2013robust} showed, using a non-rigorous analysis, that the asymptotic variance of the high-dimensional M-estimators depends on a convolution of the error and a Gaussian component. This asymptotic characterization was made rigorous by \cite{donoho2016high} using the (vector-valued) AMP framework. M-estimators in the proportional asymptotic regime have also been studied using other techniques, such as Gaussian comparison inequalities \citep{thrampoulidis2018precise} and leave-one-out analysis  \citep{lei2018asymptotics}. Finally, we mention that  AMP has been used to study the high-dimensional asymptotics of a variety of other estimators in the homoscedastic setting, including the LASSO \citep{bayati2012lasso,huang2022lasso} and logistic regression \citep{sur2019modern, zhao2022asymptotic}.

The rest of the paper is organized as follows. Section~\ref{sec: Huber mean main section} introduces the model, states the main result, and provides the corresponding numerical illustrations. Further, we demonstrate the applicability of the result by using it to construct a hypothesis test for centre  displacement. Section~\ref{sec: proof full section} develops the matrix-AMP framework and uses its state evolution to prove the main result. 
Section~\ref{sec: discussion} discusses the implications of our results and directions for future research.

\section{Model and main results}\label{sec: Huber mean main section}

Recall that the data are generated as $Y_i = X_i^\top\beta_0 + \sigma(X_i^\top \gamma_0)\varepsilon_i$, $1 \le i  \le n$, where $(X_i, \varepsilon_i)$ are independent and identically distributed copies of $(X, \bar\varepsilon)$ in \eqref{eq: data generating process}, with $X$ and $\bar\varepsilon$ independent. We consider a random design setting, with the covariates  $X_i$ and the scaled population covariance   $\Omega_p \in \reals^{p \times p}$ satisfying Assumptions~\ref{cond: covariance spectrum} and \ref{cond:white_noise_design} below. The nonlinear scale function $\sigma$ satisfies Assumption~\ref{cond: scale function}.  Throughout, limits are taken in the proportional regime where $n/p \to \delta \in (1, \infty)$ for $n, p\to\infty$.

\begin{enumerate}[label=(S{{\arabic*}})]
\item\label{cond: scale function} The scale function $\sigma:\mathbb R\to(0,\infty)$ is continuously differentiable, and for some \(\kappa>0\), the following holds for all $\alpha >0$: 
\[ E[ |\sigma(\alpha G)|^{2+\kappa} ]<\infty, \qquad E[ |\sigma'(\alpha G)|^{2+\kappa} ] <\infty , \quad \text{where } G \sim N(0,1). \]

     \item \label{cond: error distribution}
The error variable $\bar\varepsilon$ is independent of $X$, and  has zero mean and  finite second moment  $E[| \bar\varepsilon |^{2}]< \infty$. 

 \item \label{cond: covariance spectrum}
    The covariance matrices $\Omega_p$ are deterministic, positive definite, and have uniformly bounded spectrum: there exist constants $0 < c < C < \infty$, independent of $p$, such that
    \[
    c \le \lambda_{\min}(\Omega_p) \le \lambda_{\max}(\Omega_p) \le C .
    \]

    \item \label{cond:white_noise_design}
 The design matrix is $\bm X = (X_1, \ldots, X_n)^\top \in \mathbb R^{n\times p}$, 
    where the rows $X_i$ are independent and identically distributed copies of $X \in \mathbb R^p$.
    Furthermore, the rows can be expressed as $X_i = \Omega_p^{1/2} \tilde X_i $, where the whitened covariates $\tilde X_i$ have independent and identically distributed components $(\tilde X_{ij})_{1 \le j \le p}$ satisfying 
    \[
        E[\tilde X_{ij}]=0,\qquad E[\tilde X_{ij}^2]= \frac{1}{n}.
    \]
  The $\tilde{X}_{ij}$ are independent of $\bm\varepsilon = (\varepsilon_1, \ldots, \varepsilon_n)^\top$, and  $\limsup_{n \to \infty} \|\tilde X\|_{\mathrm{op}} \leq C$ almost surely,  for some constant $C$.
\end{enumerate}
Assumption~\ref{cond:white_noise_design}  is satisfied by many distributions on 
$\tilde{X}_{ij}$, including Gaussian, uniform, and Rademacher, and 
 agrees with \citet[Definition~3.1]{tan2024approximate} and \citet[Definition~2.15]{wang2024universality} for a more general class of white noise design.

 \subsection{A class of effective score functions}\label{ssec: effective score function}

Given a convex continuously differentiable loss function $\ell$ and a parameter $b \in \mathbb R_+$, let
\begin{eqnarray}\label{eq: proximal}
\bar g(u,y; \, b) = \arg\min_{z \in \mathbb R} \left\{ \ell(z, y) + \frac{1}{2 b} (z- u)^2\right\},  
\end{eqnarray}
and 
\begin{eqnarray}\label{eq: effective score}
 g(u,y; \, b) = \frac{\bar g(u,y; b) - u}{ b}  .
\end{eqnarray}
For the Huber loss in \eqref{eq: Huber loss function} with threshold $a \in (0, \infty]$, 
\[
\bar g(u, y; \, b ) = \begin{cases}
(b y + u)/(1 + b), & |y - u| \leq a (b + 1) \\
 a \cdot \operatorname{sgn} (y -u) b + u  , & |y - u| > a(b+1) ,
\end{cases}
\]
and 
\begin{eqnarray}\label{eq: g_1 function Huber}
  g(u,y; \, b) = \begin{cases}
   (y - u)/(1 + b)  , & |y - u| \leq a (b +1) \\
   a \cdot \operatorname{sgn} (y - u) ,  & |y - u| > a (b + 1).
   \end{cases}
\end{eqnarray}
The squared error loss ($a = \infty$) gives $\bar g(u,y; b) =\frac{b y + u}{1 + b}$ and
$g(u, y; \, b ) = \frac{y - u}{1 + b}$.

We will see later in \eqref{eq: asymptotic variance expression} that $g$ replaces the score function in the formula for asymptotic variance in the dense proportional regime. We will also need the (weak) partial derivatives of $g(u, y; \, b )$ with respect to $u$ and $y$. These  derivatives, denoted by  $\partial_1 g$ and $\partial_2 g$, respectively, can be readily computed from \eqref{eq: g_1 function Huber}:
\begin{align} \label{eq: g_1_derivatives}
\partial_1 g(u, y; \, b ) = - \partial_2 g(u, y; \, b ) = 
\begin{cases}
        -{1}/{(1+b)}, &   |y - u| \leq a (b + 1), \\
        0, &  |y - u| > a(b+1).
    \end{cases}
\end{align}

\subsection{Parameters for asymptotic quantification}\label{ssec: fixed point state evolution}

The parameters $(\nu_*, \tau_*^2)$ characterizing the asymptotic distribution of the Huber estimator (see \eqref{eq: asymptotic risk convergence init}) are defined as the solutions to a set of nonlinear equations, given in \eqref{eq: fixed point state evolution} below. To define these parameters, we impose the following assumption on $B = (\beta_0, \gamma_0) \in\mathbb R^{p\times 2}$.
\begin{enumerate}[resume, label=(S{{\arabic*}})]
\item \label{cond: empirical distribution of signal matrix}
 The empirical distribution of the rows of $\tilde B = \Omega_p^{1/2} B$ converges in Wasserstein-2 distance to the distribution of a random vector
    $
    \bar{B} = (\bar\beta_0,\bar\gamma_0)^\top$ with non-singular covariance matrix  and $E \|\bar{B} \|_2^2 < \infty$. 
 \end{enumerate}
Convergence in Wasserstein distance is defined in Appendix \ref{ssec: aux_results} (Definition \ref{def: wasserstein}). 

The parameters $(\nu_*,  \tau_*^2)$ are given by the solution to the following set of equations involving an additional parameter $b_* \in \reals_+$:
\begin{equation}
\label{eq: fixed point state evolution}
    \begin{split}
        \frac{1}{b_*} & =  \delta E[\partial_2 g(Z^*_1, Z_1 + \sigma(Z_2) \bar\varepsilon; \, b_*)],\\
         \nu_* & = \delta b_* \,  E[\partial_2 g (Z_1^*,  Z_1 + \sigma(Z_2) \bar\varepsilon ; \, b_* )\, \sigma'(Z_2)\bar\varepsilon], \\
        \tau_*^2 & = (\delta b_*)^2 E[g (Z_1^*,  Z_1 + \sigma(Z_2) \bar\varepsilon ; \, b_*)^2], 
    \end{split}
\end{equation}
where $(Z_1, Z_2, Z_1^*)^\top \in \reals^3$ is a zero-mean Gaussian vector independent of $\bar\varepsilon$, with covariance matrix $\bar \Sigma^*$ defined as follows. Writing $ \bar \Sigma^* = \begin{pmatrix}
    \bar\Sigma_{11} & \bar\Sigma_{12}^* \\ (\bar{\Sigma}^*_{12})^{\top} &  \bar\Sigma_{22}^*
\end{pmatrix} $, we have $\bar\Sigma_{11} = \frac{1}{\delta} E [\bar B \bar B^\top] \in \reals^{2 \times 2}$, where $\bar B = (\bar\beta_0,\bar\gamma_0)^\top$ is defined in Assumption \ref{cond: empirical distribution of signal matrix}. The other components of $\bar\Sigma^*$ are:
\begin{align}
\label{eq: barSigma*_def}
    \bar{\Sigma}^*_{12} =\frac{1}{\delta} E\Bigg[ \begin{pmatrix}
        \bar\beta_0 \\ \bar\gamma_0
    \end{pmatrix} (\bar\beta_0 + \nu_* \bar\gamma_0 + \tau_* G) \Bigg], \qquad \bar{\Sigma}^*_{22} = 
   \frac{1}{\delta} E[(\bar\beta_0 + \nu_* \bar\gamma_0 + \tau_* G)^2 ], 
\end{align}
where $G \sim N(0, 1)$ is independent of $(\bar\beta_0,\bar\gamma_0)$. 
\begin{proposition}
\label{prop: FP existence and uniqueness}
    Suppose Assumptions \ref{cond: scale function}-\ref{cond: covariance spectrum} and \ref{cond: empirical distribution of signal matrix} hold, and $\lim \frac{n}{p} := \delta >1$. Then there exists a solution $(b_*, \nu_*, \tau_*^2) \in \reals_+ \times \reals \times \reals_+$ to the system of equations in \eqref{eq: fixed point state evolution}-\eqref{eq: barSigma*_def}.
\end{proposition}
The proposition is proved in Appendix \ref{appendix: proof FP existence and uniqueness}. Moreover, if Assumption \ref{cond: restricted secant inequality} (Section~\ref{sec: proof full section}, p.\pageref{cond: restricted secant inequality}) also holds, then Theorem \ref{thm: dist convergence of huber mean} shows that the solution to \eqref{eq: fixed point state evolution}-\eqref{eq: barSigma*_def} must be unique.

The first and last equations in \eqref{eq: fixed point state evolution} are similar to the equations characterizing the asymptotic distribution of an $M$-estimator in the homoscedastic setting \citep[Eqs. (9) and (10)]{donoho2016high}. The parameter $\nu_*$ quantifies the bias in the estimator due to heteroscedasticity and, by  \eqref{eq: fixed point state evolution}, vanishes when $\sigma(\cdot)$ is a constant.

\begin{example}[Least squares]
\label{example: LS}
For the least squares estimator ($a=\infty$), we have $g(u, y; b) = (y - u)/(1 + b)$. In this case, $b_*$ can be computed in closed form and \eqref{eq: fixed point state evolution} simplifies to:
\begin{eqnarray}\label{eq: least square fixed point}
    & b_* =\frac{1}{\delta-1}, \qquad \nu_* =E[\sigma'(Z_2)\bar\varepsilon]=0, \qquad 
   & \tau_*^2 = E [(Z_1  + \sigma(Z_2)\bar\varepsilon - Z_1^*)^2] .
\end{eqnarray}
Therefore, the least-squares estimator is always correctly centred $(\nu_*=0)$ given Assumption~\ref{cond: error distribution}. In Section \ref{ssec: main result}, we give Condition~\ref{cond: a-Huber mean error} under which a general Huber estimator is correctly centred.
\end{example}

\subsection{Asymptotic quantification}\label{ssec: main result}

We now state our main result, which shows that the  Huber estimator $\hat \beta$ has the following distributional property as $p \to \infty$:  the empirical distribution of $\Omega_p^{1/2}(\hat{\beta} - \beta_0 - \nu_* \gamma_0)$ converges to the law of a Gaussian vector with distribution $N(0, \tau_*^2 I_p)$. This result is stated in terms of \emph{polynomial growth} test functions, defined in \eqref{eq: polynomial growth} (in Appendix \ref{ssec: aux_results}). Essentially,  a function $\phi: \reals^m \to \reals$ has polynomial growth of order $r$ if $|\phi(x)|$ is bounded by an $r$th-order polynomial of $\| x \|_2$ for all $x \in \reals^m$.

Since the Huber loss \eqref{eq: Huber loss function} is linear in its tails, it has vanishing second derivatives outside the quadratic region and the objective need not be strongly convex. For our analysis via AMP, we impose a \emph{restricted secant condition} anchored at the stationary AMP iterates. 
The precise statement of the condition is given in Section~\ref{sec: proof full section} (Assumption \ref{cond: restricted secant inequality} on p.\pageref{cond: restricted secant inequality}), after introducing the stationary AMP recursion. We provide a heuristic justification for this condition in Section~\ref{sec: proof full section} by establishing a lower bound on the fraction of residuals in the quadratic region near the AMP iterate. With this bound, Assumption \ref{cond: restricted secant inequality} can be interpreted as a condition on the Huber loss parameter $a$. Related approaches to handling the lack of global strong convexity include localized restricted eigenvalue conditions in sparse regression \citep{sun2020adaptive} and a vanishing $\ell_2$ regularization term used to analyse unregularized M-estimators \citep{bellec2025error}.

\begin{theorem}[Distributional convergence of the Huber mean estimator]\label{thm: dist convergence of huber mean}
Suppose Assumptions \ref{cond: scale function}-\ref{cond: empirical distribution of signal matrix} hold. Further, suppose that the empirical Huber objective satisfies Assumption~\ref{cond: restricted secant inequality} in Section~\ref{ssec: stationary Huber AMP}.  Let $\lim_{n,p\to \infty} {n}/{p} \coloneq \delta >1$, and further let 
\begin{equation} \label{eq: transformed coordinates}
    \check\beta : = \Omega_p^{1/2} \hat \beta , \qquad  \tilde\beta_{0} \coloneq  \Omega_p^{1/2} \beta_0, \qquad  \tilde\gamma_{0} \coloneq  \Omega_p^{1/2} \gamma_0.
\end{equation}
    Then, for any continuous function $\phi: \mathbb R^3 \to \mathbb R$ satisfying the polynomial growth condition with order $2$ (defined in \eqref{eq: polynomial growth}), we have almost surely
    \begin{eqnarray}\label{eq: asymptotic risk Huber estimator}
 \frac{1}{p} \sum_{j =1}^p \phi (\check\beta_j -  \tilde\beta_{0, j} - \nu_* \tilde\gamma_{0, j}, \,  \tilde\beta_{0, j}, \, \tilde\gamma_{0, j} ) {\to} E [\phi (\tau_* G, \bar\beta_0, \bar\gamma_0)],    
    \end{eqnarray}
    where $G\sim N(0,1)$ is independent of $(\bar\beta_0, \bar\gamma_0)$, and $\nu_*, \tau_*^2$ correspond to the solution of \eqref{eq: fixed point state evolution}-\eqref{eq: barSigma*_def}, which is unique. 
\end{theorem}

By taking the test function  $\phi (z, u, v)  = z^2$, Theorem \ref{thm: dist convergence of huber mean} gives the asymptotic variance.
\begin{corollary} [Asymptotic variance of the Huber estimator]\label{cor: asymptotic characterization}
Under the setting of Theorem \ref{thm: dist convergence of huber mean}, we have almost surely 
    \begin{eqnarray}\label{eq: asymptotic risk convergence}
            \frac{1}{p} (\hat\beta - \beta_0 - \nu_* \gamma_0)^\top \Omega_p (\hat\beta -  \beta_0 - \nu_* \gamma_0) {\to} \tau_*^2.   
    \end{eqnarray}
\end{corollary}

Although heteroscedasticity enters the model through the nonlinear function $\sigma(\cdot)$, we observe that its effect on the asymptotic centre  is linear in $\gamma_0$, with $\sigma(\cdot)$ determining the magnitude of the coefficient $\nu_*$. 
We now discuss a few key consequences of Theorem \ref{thm: dist convergence of huber mean} and Corollary \ref{cor: asymptotic characterization}.

Under the assumptions of Theorem \ref{thm: dist convergence of huber mean}, a sufficient condition for the Huber estimator to be \emph{Fisher consistent} for $\beta_0$ is the following:
 \begin{enumerate}[resume,label=(C1)]
    \item \label{cond: a-Huber mean error} 
   ${E} [\partial_2 g ( Z_1^*, Z_1 + \sigma(Z_2) \bar\varepsilon; b_*) \sigma'(Z_2) \bar\varepsilon ]= 0$.
\end{enumerate}

By the expression in \eqref{eq: fixed point state evolution}, Condition~\ref{cond: a-Huber mean error} implies $\nu_* = 0$ given $b_* > 0$. For the Huber loss with $\partial_2 g$ in \eqref{eq: g_1_derivatives}, Condition~\ref{cond: a-Huber mean error} can further be written as 
\[
E\left[\sigma'(Z_2)\bar\varepsilon \, 
I\big\{\left|Z_1+\sigma(Z_2)\bar\varepsilon-Z_1^*
\right|<a(1+b_*)\big\}\right]=0.
\]

Under the proposed heteroscedastic model, Condition~\ref{cond: a-Huber mean error} characterises whether the estimator remains asymptotically centred at the regression parameter or acquires a first-order component in the heteroscedastic direction. This extends the classical Fisher-consistency characterization of the Huber mean estimator to the dense proportional asymptotic regime. For least squares  ($a = \infty$), Condition~\ref{cond: a-Huber mean error} reduces to $E [\bar\varepsilon ] = 0$. For $a \in (0, \infty)$ to guarantee $\nu_* = 0$, it suffices that the error distribution is symmetric around 0.

\begin{remark}
   Conditions similar to~\ref{cond: a-Huber mean error} have been observed for the efficacy of computationally efficient estimators for phase retrieval \citep{mondelli2021approximate} and other generalized linear models \citep{damian2024computational} with random design in the proportional asymptotic regime.
\end{remark}

\emph{Asymptotic variance}. With the centre  $\beta_0+\nu_*\gamma_0$,  the asymptotic variance in Corollary \ref{cor: asymptotic characterization} can be written as follows, using \eqref{eq: fixed point state evolution}:
    \begin{eqnarray}\label{eq: asymptotic variance expression}
\tau_*^2 = \frac{E [g (Z_1^*, Z_1 +\sigma(Z_2) \bar\varepsilon)^2 ]}{(E[\partial_1 g (Z_1^*, Z_1 +\sigma(Z_2) \bar\varepsilon)])^2},    
\end{eqnarray}
where $(Z_1, Z_2, Z^*_1)^\top \sim N(0, \bar\Sigma^*)$ is independent of $\bar\varepsilon$. Comparing with the classic asymptotic variance expression $E[\psi_a^2] / E[\partial_1\psi_a]^{2} $ for $\psi_a =\partial_1\ell_a$, the Huber score function, we observe that the asymptotic variance expression in \eqref{eq: asymptotic variance expression} takes a classical sandwich form, with the ordinary Huber score $\psi_a$ replaced by the effective score function $g$. As in \citet{donoho2016high}, the effective score is evaluated under a Gaussian-perturbed error distribution. Here, however, the Gaussian perturbation is combined with the heteroscedastic error $\sigma(Z_2)\bar\varepsilon$. 
To unify with the existing literature, we use the homoscedastic least-squares model as a benchmark comparison. For this case, which  corresponds to setting $\sigma(\cdot) = 1$, $E[\bar\varepsilon] = 0$, $a = \infty$, and $\ell_a(z,y)=\frac12(y-z)^2$, the asymptotic variance expression in \eqref{eq: asymptotic variance expression}  becomes
\begin{equation}\label{eq: asymptotic variance least squares benchmark}
\tau_{*, \, {\textrm{LS}}}^2 = E\left[ \left(\bar\varepsilon +  \sqrt{{\tau_{*, \, {\textrm{LS}}}^2}/{\delta}} \, G \right)^2 \right], \quad G \sim N(0,1).
\end{equation}
This composition indicates that the asymptotic variance is a convolution of the error and an independent Gaussian component. Further, rearranging the above equation gives
\[
\tau_{*, \, {\textrm{LS}}}^2 = \frac{E[\bar\varepsilon^2]}{1- 1/\delta} .
\]
For centred Gaussian errors, it can be written as $\tau_{*, {\textrm{LS}}}^2 = \frac{1}{I(\bar\varepsilon)}\frac{1}{1- 1/\delta}$,
where $I(\bar\varepsilon)$ is the classic Fisher information for a location parameter under the error distribution $F_{\bar\varepsilon}$. In contrast with the fixed-$p$ regime, the variance does not vanish when $\delta$ is held fixed and $n,p \to \infty$. Moreover, the variance blows up as $\delta \to 1$.

\subsection{Illustrations}\label{ssec: work out example}

\begin{figure}[t!]
    \centering
    \includegraphics[width=\textwidth]{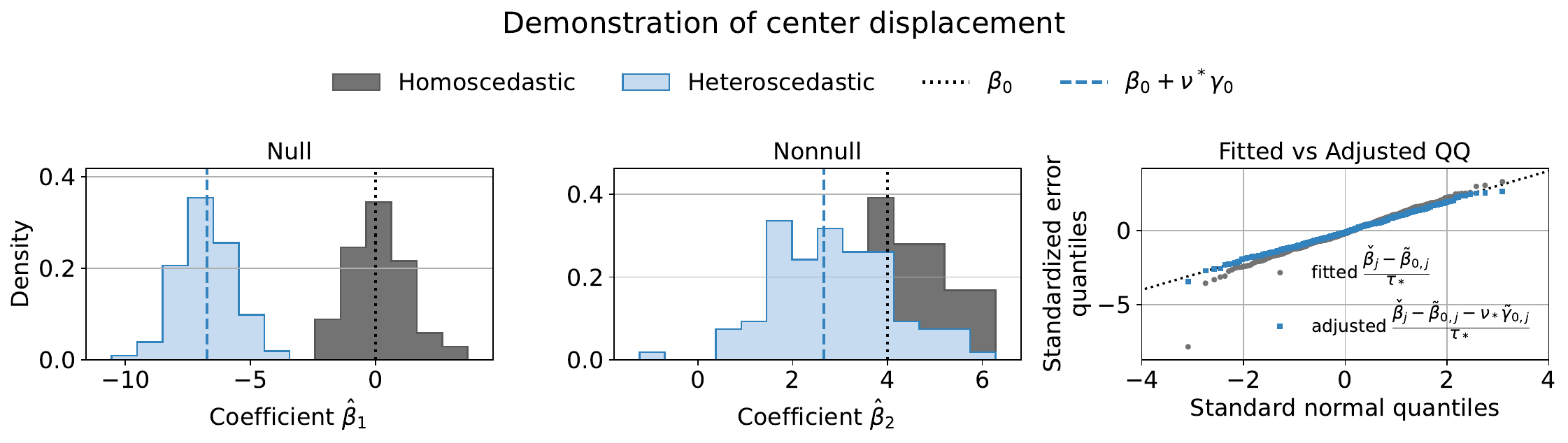}
    \caption{$Y_i = X_i^\top \beta_0 + \sigma(X_i^\top \gamma_0)\varepsilon_i$, $X_i = \Omega_p^{1/2}\tilde X_i$, $\tilde X_{ij}\sim N(0,1/n)$, and $(\Omega_p)_{jk}=0.35^{|j-k|}$. $p = 500$, $n = 4000$, $\varepsilon_i = 2 ( e_i - E(e_i))/\sqrt{\operatorname{var} (e_i)}$, where $e_i \sim  0.9 N(0, 1) + 0.1 \delta_{24}$. For heteroscedastic variance, $\sigma(z) = \exp(z)$; for homoscedastic variance, $\sigma(z) = 1$. We set $\beta_{0,1} = 0$, $\gamma_{0,1} = 10$, $\beta_{0, 2} = 4$, and $\gamma_{0,2} = 2$. The other components of $\beta_0$ and $\gamma_0$ are independently generated from $N(0,1)$. The Huber threshold is $a = 1.345$. The left and the middle panels are the histograms of $\hat\beta_{1}$ and $\hat\beta_{2}$, respectively, over 100 replications. The right panel shows QQ-plots of the fitted standardized errors $(\check{\beta}_j-\tilde{\beta}_{0,j})/\tau_*$ and the recentred standardized errors $(\check{\beta}_j-\tilde{\beta}_{0,j}-\nu_*\tilde{\gamma}_{0,j})/\tau_*$ against standard normal quantiles, using the first of the 100 simulation replications. }
    \label{fig: recenter demo} 
\end{figure}

We present numerical experiments which serve two purposes. The first is to make the centre  displacement $\nu_* \gamma_0$ visible at the coordinate level. The second is to assess the fixed-point prediction for the normalized quadratic error in \eqref{eq: asymptotic risk convergence}. All reported quantities are computed from independently generated datasets, and the deterministic parameters
$(\nu_*,\tau_*^2)$ are obtained by numerically solving \eqref{eq: fixed point state evolution}.

Figure~\ref{fig: recenter demo} isolates the displacement phenomenon in an asymmetric-error example. We use $p=500$, $n=4000$, $(\Omega_p)_{jk}=0.35^{|j-k|}$, and Huber threshold
$a=1.345$. To assess the recentring effect induced by heteroscedasticity, we compare the Huber estimator under $\sigma(z)=\exp(z)$ with the homoscedastic case ($\sigma(z)=1$), while holding all other aspects of the data-generating process fixed. For each simulation setting, the parameter vectors $\beta_0$ and $\gamma_0$ are generated once and held fixed across all replications.
The first two coordinates are fixed at
$\beta_{0,1}=0$, $\gamma_{0,1}=10$, $\beta_{0,2}=4$, and
$\gamma_{0,2}=2$, while all remaining coordinates are generated independently from $N(0,1)$. The error is a centred and standardized contaminated
distribution based on $0.9N(0,1)+0.1\delta_{24}$. The coordinate histograms show that the raw Huber estimator is displaced in the direction predicted by
$\gamma_0$. Using the first of the 100 simulation replications, the right panel shows QQ-plots of the fitted standardized errors $(\check{\beta}_j-\tilde{\beta}_{0,j})/\tau_*$ and the recentred standardized errors $(\check{\beta}_j-\tilde{\beta}_{0,j}-\nu_*\tilde{\gamma}_{0,j})/\tau_*$ against standard normal quantiles. We observe that subtracting $\nu_*\gamma_0$ removes the centre  displacement.

To illustrate \eqref{eq: asymptotic risk convergence} for the variance comparison, we evaluate the following quantities: 
\begin{eqnarray}\label{eq: MSE ratios}
    & R_{\rm centred} =\left\{\frac{1}{p} (\hat\beta -  \beta_0 - \nu_* \gamma_0)^\top \Omega_p (\hat\beta -  \beta_0 - \nu_* \gamma_0)\right\} / \tau_*^2 ,  \nonumber\\
    & R_{\rm uncentred} =\left\{\frac{1}{p} (\hat\beta -  \beta_0 )^\top \Omega_p (\hat\beta -  \beta_0)\right\} / \tau_*^2,
\end{eqnarray}
where $\hat\beta$ is the Huber estimator obtained by solving \eqref{eq: M-estimator minimization}. By Corollary~\ref{cor: asymptotic characterization}, $R_{\rm centred}$, which includes the centre  displacement adjustment,  should fluctuate around 1, and $R_{\rm uncentred}$  reveals whether a centre  displacement correction is needed.

\begin{figure}[t!]
    \centering
    \includegraphics[width=0.85\textwidth]{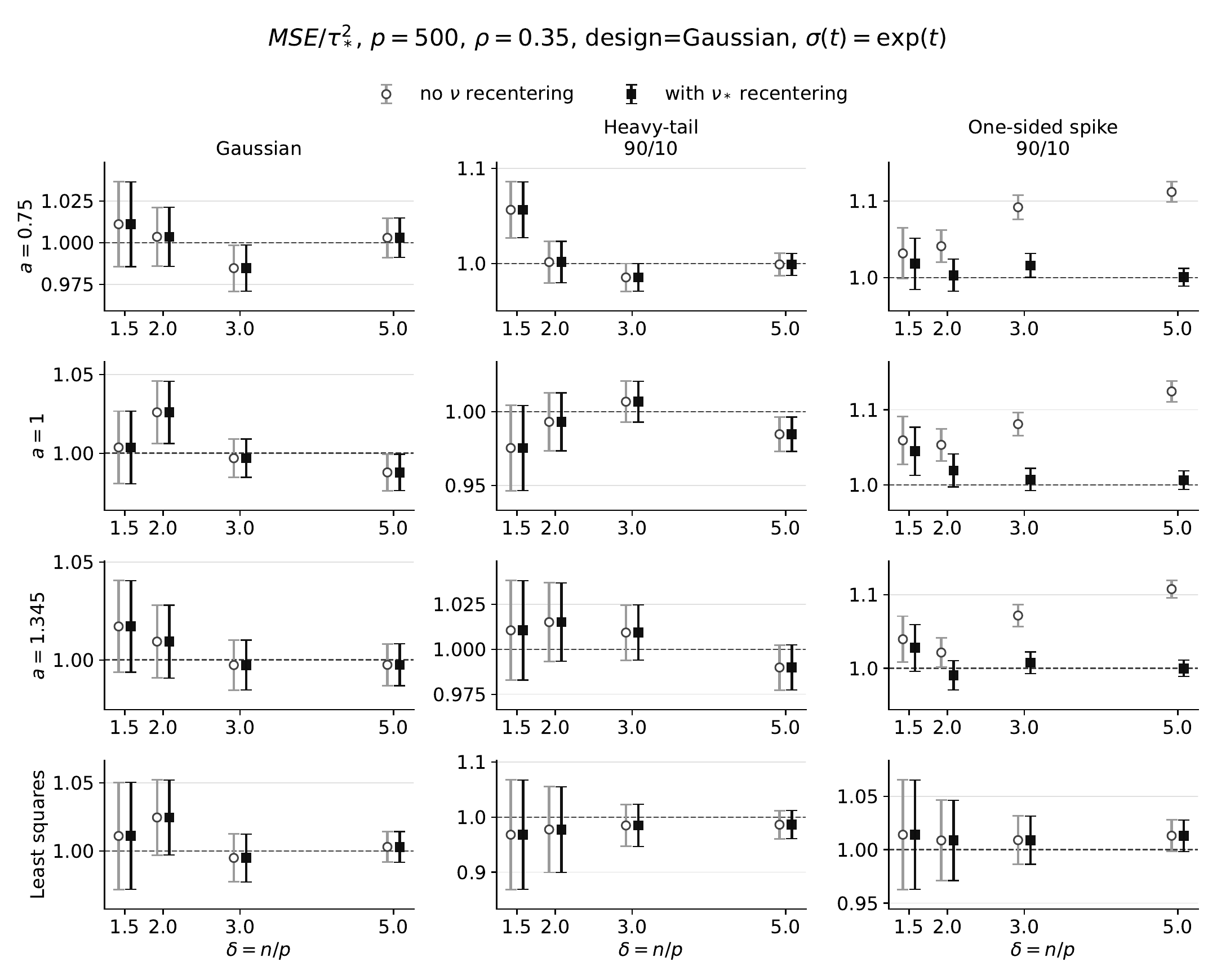}
    \caption{$Y_i = X_i^\top \beta_0 + \exp(X_i^\top \gamma_0)\varepsilon_i$, $X_i = \Omega_p^{1/2}\tilde X_i$, $\tilde X_{ij}\sim N(0,1/n)$, and $(\Omega_p)_{jk}=0.35^{|j-k|}$. $p = 500$, $\delta = 1.5, 2, 3, 5$. The error $\varepsilon_i =  ( e_i - E(e_i))/\sqrt{\operatorname{var} (e_i)}$, where $e_i \sim N(0,1)$ (left column), $0.9N(0, 1) + 0.1 N(0, 16)$ (middle column), $0.9 N(0, 1) + 0.1 \delta_{8}$ (right column). The parameter vectors $\beta_0$ and $\gamma_0$ are generated once and held fixed across all replications. We set $\beta_{0,1} = 0$, $\gamma_{0,1} = 10$, $\beta_{0, 2} = 4$, and $\gamma_{0,2} = 2$. The other components of $\beta_0$ and $\gamma_0$ are randomly generated from a standard Gaussian distribution $N(0,1)$. The rows correspond to Huber threshold $a = 0.75, 1, 1.345, \infty$, respectively. The grey and black error bars (mean and 1.96 standard deviations over 200 replications) represent the ratios $R_{\rm uncentred}$ and $R_{\rm centred}$ in \eqref{eq: MSE ratios}. }
    \label{fig: gaussian exp}
\end{figure}

\begin{figure}[t!]
    \centering
    \includegraphics[width=0.85\textwidth]{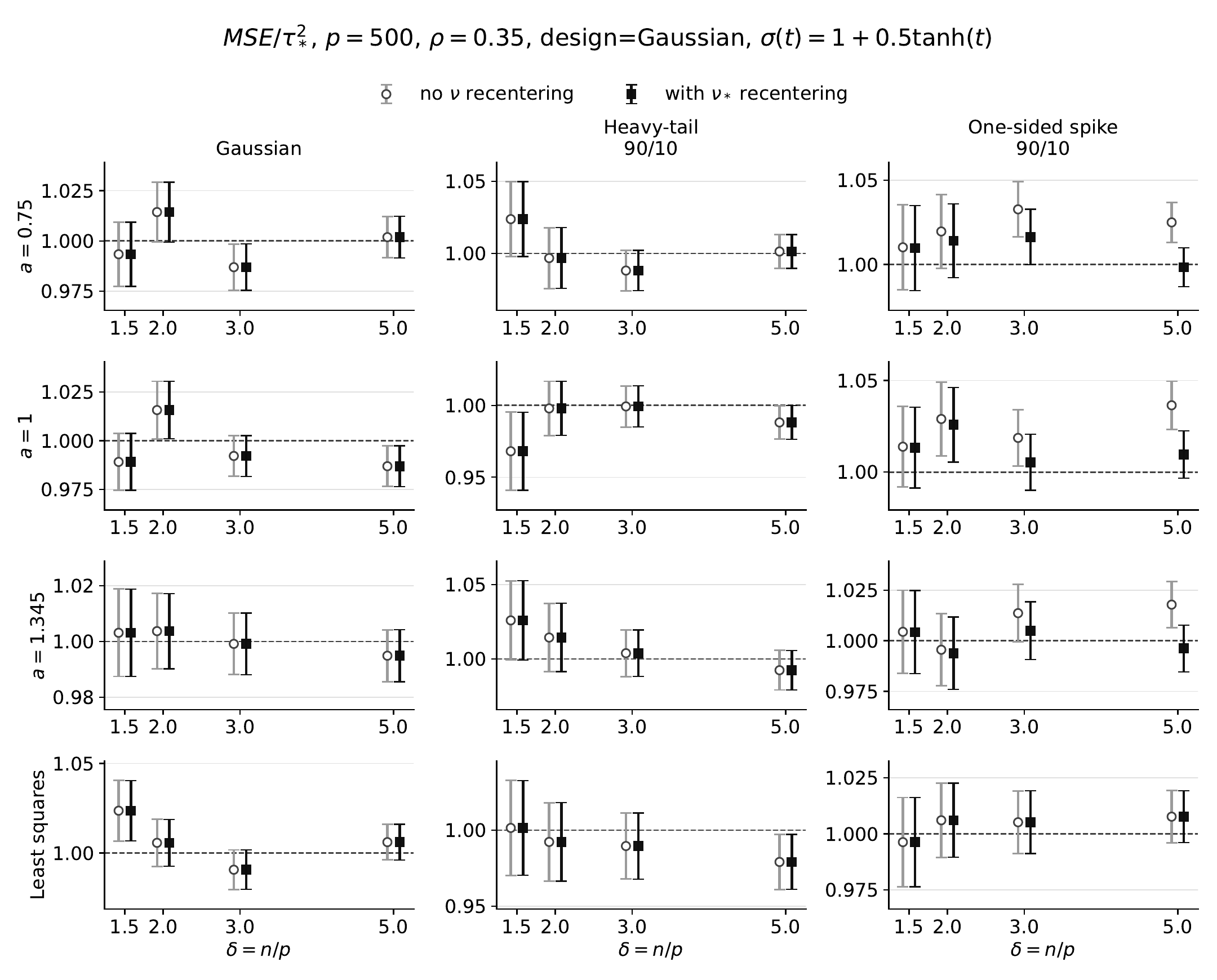}
    \caption{$Y_i = X_i^\top \beta_0 + (1 + 0.5\tanh(X_i^\top \gamma_0))\varepsilon_i$, $X_i = \Omega_p^{1/2}\tilde X_i$, $\tilde X_{ij}\sim N(0,1/n)$, and $(\Omega_p)_{jk}=0.35^{|j-k|}$. $p = 500$, $\delta = 1.5, 2, 3, 5$. The error $\varepsilon_i =  ( e_i - E(e_i))/\sqrt{\operatorname{var} (e_i)}$, where $e_i \sim N(0,1)$ (left column), $0.9N(0, 1) + 0.1 N(0, 16)$ (middle column), $0.9 N(0, 1) + 0.1 \delta_{8}$ (right column). The parameter vectors $\beta_0$ and $\gamma_0$ are generated once and held fixed across all replications. We set $\beta_{0,1} = 0$, $\gamma_{0,1} = 10$, $\beta_{0, 2} = 4$, and $\gamma_{0,2} = 2$. The other components of $\beta_0$ and $\gamma_0$ are randomly generated from a standard Gaussian distribution $N(0,1)$. The rows correspond to Huber threshold $a = 0.75, 1, 1.345, \infty$, respectively. The grey and black error bars (mean and 1.96 standard deviations over 200 replications) represent the ratios $R_{\rm uncentred}$ and $R_{\rm centred}$ in \eqref{eq: MSE ratios}. }
    \label{fig: gaussian tanh}
\end{figure}

In Figures~\ref{fig: gaussian exp} and \ref{fig: gaussian tanh}, we examine the above ratios for different error distributions and varying $\delta=\frac{n}{p}$.   We generate $K=200$ independent datasets with $p=500$ and $\delta \in\{1.5,2.0,3.0,5.0\}$, and fit the corresponding Huber estimators. In all simulations, the scaled population covariance is specified by
$(\Omega_p)_{jj'}=0.35^{|j-j'|}, j,j'\in[p].
$
In Figure~\ref{fig: gaussian exp}, the variance function is $\sigma(z)=\exp(z)$,  and in Figure~\ref{fig: gaussian tanh}  we use $\sigma(z)=1+0.5\tanh(z)$.  We consider three error distributions: $N(0,1)$, $0.9N(0,1)+0.1N(0,16)$ and $0.9N(0,1)+0.1\delta_8$, to demonstrate when Condition~\ref{cond: a-Huber mean error} holds and when it is violated. 

By Corollary~\ref{cor: asymptotic characterization}, the ratio $R_{\rm centred}$ is expected to fluctuate around one for all error distributions. This is confirmed by Figures~\ref{fig: gaussian exp} and \ref{fig: gaussian tanh}. For symmetric error distributions, Condition~\ref{cond: a-Huber mean error} holds, so $\nu_* = 0$. That is, the Huber estimator is centred at $\beta_0$ and $R_{\rm uncentred}$ is also expected to fluctuate around one. This is seen in the first two columns of Figures~\ref{fig: gaussian exp} and \ref{fig: gaussian tanh}, where the error bars coincide.  For asymmetric error distributions, where Condition~\ref{cond: a-Huber mean error} is violated, the Huber estimator exhibits a centre displacement $\nu_*\gamma_0$.   The final column shows that this displacement effect is accurately captured by the coefficient $\nu_*$. After recentring by $\beta_0+\nu_*\gamma_0$, the normalized ratio again fluctuates around one, in agreement with the theoretical prediction. We observe that the least-squares estimator is centred for any error distribution, consistent with \eqref{eq: least square fixed point} which shows that $\nu_*=0$.

Throughout this section, the design matrix is Gaussian. Simulations with uniform and Rademacher designs yield qualitatively identical conclusions and are reported in Appendix~\ref{appendix: figures}. This is consistent with Theorem \ref{thm: dist convergence of huber mean}, which holds for any whitened covariate distribution satisfying Assumption \ref{cond:white_noise_design}.

\subsection{Hypothesis testing for centre displacement $\nu_*\gamma_0$}
A natural application of Theorem~\ref{thm: dist convergence of huber mean} is to test for centre displacement in the Huber estimator. That is, we test between the hypotheses $H_0: \nu_{*}\gamma_0 = 0$ versus $H_1: \nu_{*} \gamma_0\neq 0 $.

Let $\Omega_p$ satisfy Assumption~\ref{cond: covariance spectrum}. We consider a Huber estimator $\hat\beta_a$ for fixed $a\in (0, \infty)$ and a least-squares estimator $\hat\beta_{\rm LS}$, noting that $\nu_*=0$ for least-squares (Example \ref{example: LS}). Recalling that $\tilde\beta_0=\Omega_p^{1/2}\beta_0$ and $\tilde\gamma_0=\Omega_p^{1/2}\gamma_0$, by Theorem~\ref{thm: dist convergence of huber mean}, 
\begin{eqnarray}\label{eq: heuristic Gaussian approximation}
   \check\beta_{a,j} =  (\Omega_p^{1/2}\hat\beta_a)_{j} \stackrel{d}{\approx} N (\tilde\beta_{0,j} + \nu_* \tilde \gamma_{0,j}, \tau_{a , *}^2 ) , \qquad
\check\beta_{{\rm LS}, j} = (\Omega_p^{1/2}\hat\beta_{\rm LS})_j \stackrel{d}{\approx} N (\tilde \beta_{0,j} , \tau_{\rm LS, *}^2 ). 
\end{eqnarray}

Let $D=\hat\beta_a-\hat\beta_{\rm LS}=\Omega_p^{-1/2}(\check\beta_a-\check\beta_{\rm LS})$. The Gaussian approximations in \eqref{eq: heuristic Gaussian approximation} suggest that the mean shift of $D$ is $\nu_*\gamma_0$. The covariance of $D$ must account for the correlation between $\hat\beta_a$ and $\hat\beta_{\rm LS}$ and for the unknown $\Omega_p$. Since $\|\nu_*\gamma_0\|_2^2=\nu_*^2\gamma_0^\top\gamma_0$,
the null $H_0$ is equivalent to the vanishing of this shift when $\gamma_0^\top\gamma_0>0$, motivating a statistic based on $\|D\|_2^2$.
If $\gamma_0=0$, the shift vanishes irrespective of $\nu_*$.

A statistic based on $\|D\|_2^2$ requires a correction for the noise contribution due to the correlation between $\hat\beta_{a}$ and $\hat\beta_{\rm LS}$. To avoid estimating this contribution for centring, we consider a normalized split-sample cross product motivated by the sample-splitting principle underlying \citet{kim2024dimension}. 
We use this construction to illustrate the applicability of our main result; procedures that use the data more efficiently may offer further gains in power. We split the observations into two disjoint sets
$I_l=\{(l-1) \lfloor n/2\rfloor +1,\ldots,l \lfloor n/2\rfloor\},l=1,2$.
We fit both estimators on each subset, denoting them by $\hat\beta_a^{(l)}$ and $\hat\beta_{\rm LS}^{(l)}$, and define $D_l=\hat\beta_a^{(l)}-\hat\beta_{\rm LS}^{(l)}$, $l=1,2$. Conditional on the coefficient vectors, $D_1$ and $D_2$ are independent. Heuristically, the unscaled cross-product $D_1^\top D_2$ estimates $\nu_{*,h}^2\gamma_0^\top\gamma_0$, where $\nu_{*,h}$ is the nuisance coefficient for the half sample with ratio $\delta/2$. Here, we assume $\delta >2$ to ensure that the aspect ratio of each half-sample exceeds one. 
Although the values of $\nu_*$ and $\nu_{*,h}$ may differ because they correspond to different aspect ratios, $\nu_*=0$ if and only if $\nu_{*,h}=0$. 
We therefore consider the half-sample hypotheses
\[ H_0^{(h)}: \nu_{*,h}\gamma_0 = 0 \qquad \mbox{versus}\qquad H_1^{(h)}: \nu_{*,h} \gamma_0\neq 0 .
 \]

Under the half-sample null, we have $\nu_{*,h}=0$ and the Gaussian approximation,
$E[D_l]=0$ and $E(D_l D_l^\top)=\xi_{*,h}\Omega_p^{-1}$, where $\xi_{*,h}>0$ is the covariance scale of the difference. Independence of $D_1$ and $D_2$ then gives
\[
\operatorname{var}(D_1^\top D_2)
=\xi_{*,h}^2\operatorname{tr}(\Omega_p^{-2}),
\qquad
E[\|D_l\|_2^2]
=\xi_{*,h}\operatorname{tr}(\Omega_p^{-1}).
\]
Under Assumption~\ref{cond: covariance spectrum}, the squared norms concentrate around their expectations. Consequently, the asymptotic variance factor for the normalized cross-product is
\[
r=
\frac{\operatorname{tr}(\Omega_p^{-2})}
{\{\operatorname{tr}(\Omega_p^{-1})\}^2},
\]
in which the unknown scale $\xi_{*,h}$ cancels. Intuitively, under the null, $\|D_1\|_2$ and $\|D_2\|_2$ account for the overall noise levels of the two vectors, while $r$ adjusts for the covariance structure induced by the unknown $\Omega_p$.

Letting $W=n X^\top X$, the inverse-Wishart identities of \citet{haff1979identity} motivate estimating $\operatorname{tr}(\Omega_p^{-1})$ by $(n-p-1)\operatorname{tr}(W^{-1})$ and $\operatorname{tr}(\Omega_p^{-2})$ by $(n-p)(n-p-3)\operatorname{tr}(W^{-2})
-(n-p-1)\{\operatorname{tr}(W^{-1})\}^2$, which leads to a simplified estimator of $r$ as follows:
\[\hat r =
\frac{\operatorname{tr}(W^{-2})}
{\{\operatorname{tr}(W^{-1})\}^2} -\frac{1}{n-p}.\]

Thus, the final test statistic is defined as
\begin{eqnarray*}
    T =  \hat r^{-1/2}\frac{ D_1^\top D_2}
{\sqrt{D_1^\top D_1}\sqrt{D_2^\top D_2}},
\end{eqnarray*}
which we expect to converge to $N(0,1)$ in distribution. Thus the null $H_0$ is rejected at level $\alpha$ if $T > \Phi^{-1} (1- \alpha)$ and is interpreted as heteroscedasticity along a variance index parameter such that Condition~\ref{cond: a-Huber mean error} fails. The Huber estimator should then be interpreted with care, as its asymptotic centring may differ from $\beta_0$.

\begin{figure}[t!]
    \centering
    \includegraphics[width=\textwidth]{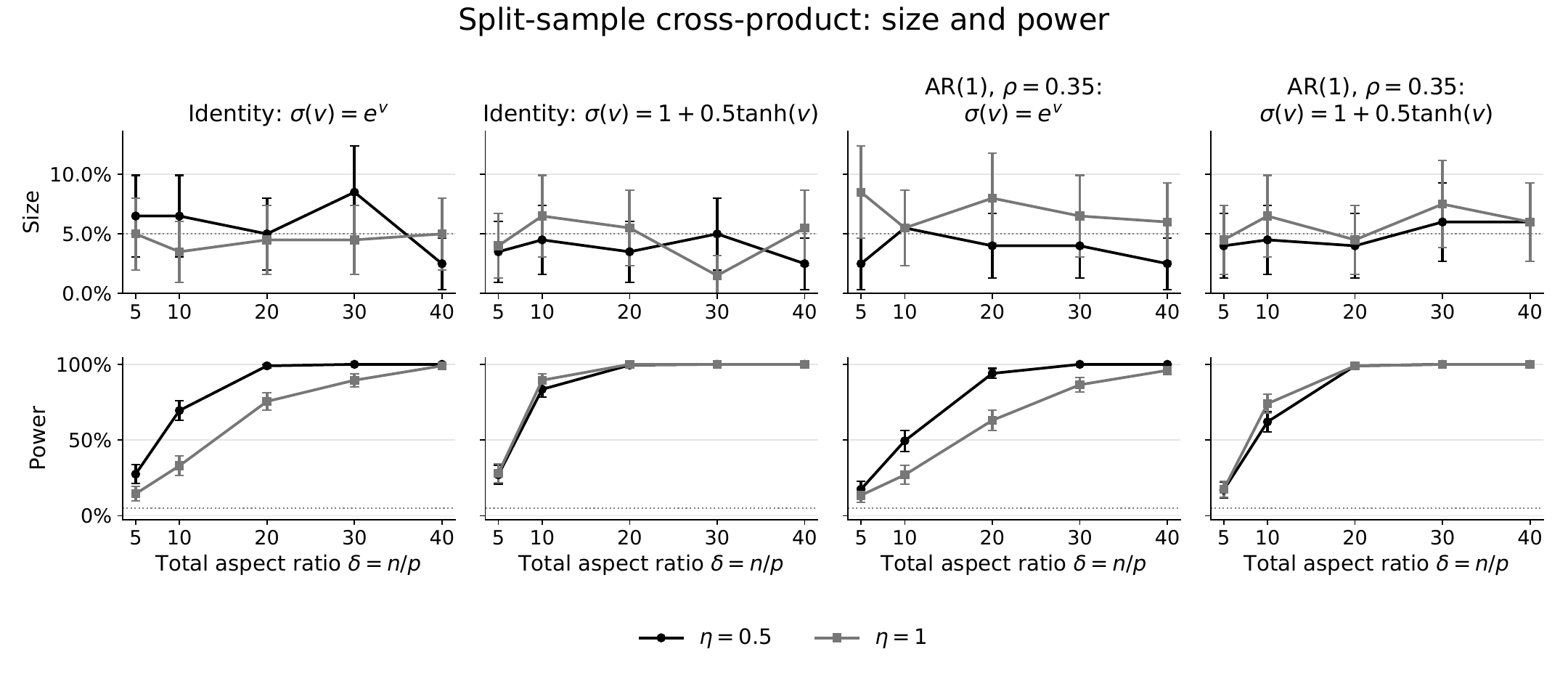}
    \caption{$Y_i = X_i^\top \beta_0 + \sigma(X_i^\top \gamma_0)\varepsilon_i$, $X_i = \Omega_p^{1/2}\tilde X_i$, $\tilde X_{ij}\sim N(0,1/n)$.  Covariance:  $\Omega_p = I_p$ (columns 1\&2) and $(\Omega_p)_{jk}=0.35^{|j-k|}$ (columns 3\&4). $p = 500$, $\delta = 5 , 10, 20, 30 ,40$. Scale function: $\sigma(z) = \exp (z)$ (columns 1\&3) and $\sigma(z) = 1 + 0.5 \tanh (z)$ (columns 2\&4). The error $\varepsilon_i =  ( e_i - E(e_i))/\sqrt{\operatorname{var} (e_i)}$, where $e_i \sim 0.9N(0, 1) + 0.1 N(0, 16)$ (top row), and $e_i \sim 0.9 N(0, 1) + 0.1 \delta_{8}$ (bottom row). The regression and variance index coefficients are drawn independently once for each simulation setting and held fixed across all replications. Specifically, $\beta_{0,j}\sim N(0,\delta)$ and $\gamma_{0,j}\sim N(0,\delta\eta)$, with $\eta\in\{0.5,1\}$. In each replication, errors are centred and scaled to unit variance.
    The Huber threshold is $a=1.345$. Points show rejection proportions over 200 replications; bars indicate 1.96 standard errors. }
    \label{fig: test statistic}
\end{figure}

The simulation settings are similar to those in Section~\ref{ssec: work out example}. We take $p=500$ and $\delta\in\{5,10,20,30,40\}$. For each setting, we independently generate $\gamma_{0,j}$ from $N(0,\delta\eta)$ for $\eta\in\{0.5,1\}$ and $\beta_{0,j}$ from $N(0,\delta)$, and fix them across replications. In each replication, the error $\varepsilon_i =  ( e_i - E(e_i))/\sqrt{\operatorname{var} (e_i)}$, where the symmetric $0.9N(0, 1) + 0.1 N(0, 16)$ and the asymmetric $0.9 N(0, 1) + 0.1 \delta_{8}$ error distributions are considered. We take $X_i = \Omega_p^{1/2}\tilde X_i$, $\tilde X_{ij}\sim N(0,1/n)$. We consider both $\Omega_p = I_p$ and $(\Omega_p)_{jk}=0.35^{|j-k|}$. Results in Figure~\ref{fig: test statistic} are based on $K=200$ replications and the points show the mean rejection proportions, with error bars indicating $1.96$ standard errors on either side. The Huber estimator uses $a=1.345$. Recall that, under Condition~\ref{cond: a-Huber mean error}, a symmetric error distribution ensures the Fisher consistency of the Huber mean estimator and a vanishing nuisance coefficient $\nu_* = 0$. The top row, based on symmetric error, shows the empirical size of the test, while the bottom row, based on asymmetric error, demonstrates its empirical power. The size stays close to the nominal level $\alpha = 0.05$. For bounded $\sigma(z) = 1 + 0.5 \tanh (z)$, the power increases for stronger signal strength ($\eta = 1$ vs $\eta = 0.5$). The opposite pattern occurs for $\sigma(t)=\exp(t)$, possibly because noise variance grows faster than signal variance, reducing the effective signal-to-noise ratio. 

\section{Proof of Theorem \ref{thm: dist convergence of huber mean}}\label{sec: proof full section}

The main idea in the proof is to establish an equivalence between the Huber estimator and a certain matrix-valued Approximate Message Passing (AMP) iteration. Working with whitened coordinates, the proof proceeds as follows. We first introduce a generic matrix AMP algorithm whose iterates, which can be interpreted as estimates of $(\beta_0, \gamma_0)$, have a succinct distributional characterization in the proportional asymptotic limit. This characterization is in terms of a set of parameters specified by a deterministic `state evolution' recursion. We then specialize the matrix AMP algorithm so that its iterates track the effect of $\beta_0$ and $\gamma_0$ on the Huber estimator. With a carefully chosen initialization, we show that the AMP iterates converge to the Huber estimator. This allows us to translate the AMP state-evolution guarantees to the Huber estimator, yielding Theorem \ref{thm: dist convergence of huber mean}.

This proof strategy is similar to other AMP-based analyses of high-dimensional estimators, e.g.,  \cite{donoho2016high, sur2019modern}, with the key difference that we use a matrix-valued AMP  \citep{tan2024approximate} to track the effect of both $\beta_0$ and $\gamma_0$ on the Huber estimator. We provide the necessary background in Section \ref{ssec: generic matrix AMP}. We refer the reader to \citet{feng2022unifying} for a detailed review of AMP, including a description of how it can be used to study the high-dimensional asymptotics of convex estimators for linear and generalized linear models.

  \subsection{Whitening transformation}\label{ssec: whitening}

Recall from Assumption~\ref{cond:white_noise_design} that $X_i = \Omega_p^{1/2}\tilde X_i$,  where $\tilde X_i$ is the whitened covariate vector. The heteroscedastic regression model can be reparameterized as
\begin{equation}\label{eq: standardized data generating process}
    Y_i=\tilde X_i^\top\tilde\beta_0
    +\sigma(\tilde X_i^\top\tilde\gamma_0)\varepsilon_i,
\end{equation}
where $\tilde \beta_0, \tilde\gamma_0$ are the pseudo-parameters defined in \eqref{eq: transformed coordinates} that will be used throughout the proof. We write
$B = (\beta_0, \gamma_0) \in\mathbb R^{p\times 2}, 
\tilde B=(\tilde\beta_0,\tilde\gamma_0)\in\mathbb R^{p\times 2}$, and $\Theta_i=(\theta_{1i},\theta_{2i})^\top= B^\top X_i = \tilde B^\top\tilde X_i$.
Then, the  model in \eqref{eq: standardized data generating process} can further be expressed as  $Y_i=\theta_{1i}+\sigma(\theta_{2i})\varepsilon_i$.

In the whitened coordinates, the regression parameter $\tilde\beta_0$ can be estimated by solving 
\begin{eqnarray}\label{eq: minimization whitened coordinate}
      \check{\beta}\in\arg\min_{ \beta \in\mathbb R^p}
 L_n (\beta), \quad  L_n (\beta)\coloneq
    \sum_{i=1}^n\ell_a(\tilde X_i^\top \beta ,Y_i).
\end{eqnarray}
Comparing \eqref{eq: minimization whitened coordinate} and \eqref{eq: M-estimator minimization}, the estimators in the whitened and original coordinates are related as $\check{\beta}=\Omega_p^{1/2}\hat\beta$. The displacement effect and the asymptotic variance of the $a$-Huber mean estimator are analysed using whitened coordinates and then translated back to the original coordinates.

\subsection{Generic matrix AMP iteration}\label{ssec: generic matrix AMP}

Let $\bm Y = (Y_1, \ldots, Y_n )^\top \in \reals^n$, $\bm X = (X_1, \ldots, X_n)^\top \in \mathbb R^{n\times p}$, and $\tilde{\bm X} = \bm X\Omega_p^{-1/2}$. Moreover, let $\tilde \Theta = \tilde{\bm X} \tilde B$, where we recall from Assumption \ref{cond: empirical distribution of signal matrix} that $\tilde{B} = (\tilde\beta_0, \, \tilde\gamma_0)  \in \reals^{p \times 2}$.

With an initialization $\check B^0 = (\check\beta^0, \check\gamma^0) \in \mathbb R^{p\times 2}$ independent of $\tilde{\bm X}$ and $\widehat R^{-1} = 0$, the AMP iterates are defined as follows, for $k \ge 0$: 
\begin{eqnarray}\label{eq: generic AMP algorithm}
  &\Theta^k = \tilde{\bm X} \check B^k - \widehat R^{k-1} (\bar F^k)^\top,\ \ \ \ \  \widehat R^k  =  g_k (\Theta^k, \bm Y) \\ \nonumber
  & B^{k+1} = \tilde{\bm X}^\top \widehat R^k - \check B^k (\bar C^k)^\top, \ \ \ \ \ \check B^{k+1} = f_{k+1} (B^{k+1}),
\end{eqnarray}
where $f_k: \mathbb R^2 \to \mathbb R^2$ and $g_k: \mathbb R^2 \times \mathbb R \to \mathbb R^2$ are  user-defined Lipschitz functions that act row-wise on matrices.  The deterministic matrices $\bar F^k, \bar C^k \in \reals^{2 \times 2}$ are defined later in terms of $f_k, g_k$, respectively.  In Section~\ref{ssec: matrix AMP Huber}, we will specify the choices of $f_k$ and 
$g_k$ to study the Huber mean estimator.  For each $k \ge 1$, the matrix AMP algorithm produces iterates $B^k$ and $\Theta^k$ that may be interpreted as estimates of $\tilde B$ and $\tilde \Theta$, respectively.

The `memory' terms $- \widehat R^{k-1} (\bar F^k)^\top$ and $- \check B^k (\bar C^k)^\top$ debias the iterates, ensuring that $ \Theta^k, B^{k+1}$ each have succinct asymptotic distributional characterizations. Specifically,  from \citet[Theorem 3.4]{tan2024approximate},  we have that for each $k \ge 1$, the row-wise joint empirical distribution of $(\tilde B, B^k) \in  \reals^{p \times 4}$ converges to the joint law of the random vector $ \begin{pmatrix} \bar B \\ M_B^k \bar B + G _B^k \end{pmatrix} \in  \reals^{4}$. Here, $G_B^k \sim N(0, T_B^k)$ is independent of $\bar B =(\bar \beta_0, \bar\gamma_0) \in \reals^2$ (defined in Assumption \ref{cond: empirical distribution of signal matrix}),  with the deterministic  mean and covariance matrices  $M^k_B, T_B^k \in \reals^{2 \times 2}$  defined  via the state evolution recursion below. Similarly, the row-wise joint empirical distribution of $(\tTheta, \Theta^k) \in \reals^{n \times 4}$ converges to the law of $\begin{pmatrix} Z \\ Z^k \end{pmatrix} \sim N(0, \bar\Sigma^k)$, where $\bar\Sigma^k \in \reals^{4 \times 4}$ is defined via state evolution.

\emph{State Evolution}.  For $s=(s_1,s_2)^\top \in \reals^2$, let
\begin{equation}\label{eq: general DGP}
    q(s,v)=s_1+\sigma(s_2)v,  
\end{equation}
and define the function $h_k: \reals^2 \times \reals^2 \times \reals \to \reals^2 $ as $h_k(s, w, v) =  g_k (w, q(s, v))$, where $g_k$ is the function used in the AMP iteration \eqref{eq: generic AMP algorithm}. 

State evolution iteratively computes the matrices $M_B^k, T_B^k, \bar{\Sigma}^k$ for $k\ge 1$, starting from an initialization $\bar{\Sigma}^0$.  The initialization  $\bar{\Sigma}^0$ is given in Assumption~\ref{cond: existence initialization}, which states that the AMP initialization $\check B^0 \in \reals^{p \times 2}$ satisfies 
 \begin{eqnarray*}
       \frac{1}{n} \begin{pmatrix}
       \tilde B^\top \tilde B & \tilde B^\top  \check B^0\\
       (\check B^0)^\top \tilde B & (\check B^0)^\top (\check B^0)
       \end{pmatrix} \stackrel{a.s.}{\to} \bar{\Sigma}^0. 
\end{eqnarray*}
Then, for $k\ge 0$, with $\begin{pmatrix} Z \\ Z^k \end{pmatrix} \sim N(0, \bar{\Sigma}^k)$,  state evolution proceeds as:
\begin{eqnarray}\label{eq: state-evolution matrix AMP}
    &&M_B^{k+1} = E[D_1 h_k(Z, Z^k, \bar\varepsilon)] ,  \qquad 
    T_B^{k+1} = E[h_k(Z, Z^k, \bar\varepsilon) h_k(Z, Z^k, \bar\varepsilon)^\top],  \nonumber\\
    && \bar{\Sigma}^{k+1} = \begin{pmatrix}
        \bar{\Sigma}_{11} & \bar{\Sigma}_{12}^{k+1} \\
        \bar{\Sigma}_{21}^{k+1} & \bar{\Sigma}_{22}^{k+1}
    \end{pmatrix},
\end{eqnarray}
where $\bar{\Sigma}_{11} = \frac{1}{\delta} E[\bar{B} \bar{B}^\top]$ is the signal covariance in whitened coordinates, and 
\begin{equation}
\label{eq:Sigmak_update}
\begin{split}
      & \bar{\Sigma}_{12}^{k+1} = (\bar{\Sigma}_{21}^{k+1})^\top =  \frac{1}{\delta} E[\bar{B} f_{k+1}(M_B^{k+1} \bar{B} + G_B^{k+1})^\top], \\
    & \bar{\Sigma}_{22}^{k+1} = \frac{1}{\delta} E[ f_{k+1}(M_B^{k+1} \bar{B} + G_B^{k+1}) (f_{k+1}(M_B^{k+1} \bar{B} + G_B^{k+1}))^\top].
\end{split}
\end{equation}
In \eqref{eq:Sigmak_update}, $\bar B$ is independent of $G_B^{k+1}\sim N(0, T_B^{k+1})$.

With these definitions, the matrices $\bar{F}^k$ and $\bar{C}^k$ are defined as:
\begin{eqnarray}\label{eq: CF limiting}
\bar F^k = \frac{1}{\delta} E[D  f_k(M_B^{k} \bar{B} + G_B^{k})], \qquad \bar C^k = E [D_2h_k (Z, Z^k, \bar\varepsilon)],  
\end{eqnarray} 
where  $D f_{k}$ denotes the Jacobian of $f_{k}$, and $D_2 h_k$ denotes the Jacobian with respect to the second argument. 

The Wasserstein convergence result for $B^{k}$ and $\Theta^k$, which follows from \citet[Theorem 3.4]{tan2024approximate}, is stated in Lemma \ref{lem: wasserstein convergence of matrix AMP}.

\subsection{Specialization to Huber Regression}\label{ssec: matrix AMP Huber}

We now choose  the functions $f_k$ and $g_k$ in \eqref{eq: generic AMP algorithm} so that the fixed point of the AMP iteration  coincides with the Huber estimator. With scalars $b_k >0$  to be defined later, for $k \ge 0$, let
\begin{equation} \label{eq: gk_fk_Huber}
  g_k(w, y) = \begin{pmatrix}
         g(w_1, y ; \, b_k) \\ 0
    \end{pmatrix}, \qquad    
    f_{k+1}(w) = \begin{pmatrix}
         (\delta b_k) w_1 \\ 0
    \end{pmatrix}, \qquad w =(w_1, w_2)^\top \in \reals^2,
\end{equation}
where $g(\cdot, \cdot ; \, b_k)$ is defined in \eqref{eq: g_1 function Huber}.  With these functions, the AMP iteration in \eqref{eq: generic AMP algorithm}, which we call the `Huber AMP', has iterates of the form $\Theta^k= (\Theta^k_1, 0)$ and $B^{k+1}= (B^{k+1}_1, 0)$.  The second column of these iterates being zero reflects the fact that the Huber estimator estimates only the first column of $B = (\beta_0, \gamma_0)$. The Huber AMP updates for $\Theta^k_1 \in  \reals^n$ and $B_1^{k+1} \in \reals^p$ are:
\begin{equation}
\label{eq: Huber_AMP}
    \begin{split}
        \Theta^k_1 &= (\delta b_{k-1}) \tbX B^k_1  \, \, - b_{k-1} g(\Theta^{k-1}_1, \bm Y ; \, b_{k-1}), \\
        B^{k+1}_1 & = \tbX^\top g(\Theta^k_1, \bm Y ; \, b_k) \,  -  \,  (\delta b_{k-1}) E[\partial_1 g(Z_1^k, \,  Z_1 + \sigma(Z_2) \bar\varepsilon; \, b_k)] B^k_1.
    \end{split}
\end{equation}
And the state evolution equations in \eqref{eq: state-evolution matrix AMP} become:
\begin{equation}
\label{eq:MBTB_k1}
    \begin{split}
        M_B^{k+1} &  = \begin{pmatrix}
            E[ \partial_2 g(Z_1^k,   Z_1 + \sigma(Z_2) \bar\varepsilon; \, b_k) ] & 
            E[ \partial_2 g(Z_1^k,   Z_1 + \sigma(Z_2) \bar\varepsilon; \, b_k) \sigma'(Z_2) \bar\varepsilon] \\ 
            0 & 0 
        \end{pmatrix}, \\
        T_B^{k+1} & = \begin{pmatrix}
            E[g(Z_1^k,   Z_1 + \sigma(Z_2) \bar\varepsilon; \, b_k)^2] & 0 \\
            0 & 0 
        \end{pmatrix}.
    \end{split}
\end{equation}
Noting that the second component $Z^k_2=0$, we denote $(Z_1, Z_2, Z^k_1)^\top \sim N(0, \bar\Sigma^k)$, where the submatrices of $\bar \Sigma^{k+1} \in \reals^{3 \times 3}$ are $\bar\Sigma_{11} = \frac{1}{\delta} E[\bar{B} \bar{B}^\top]$ and 
\begin{equation}
\label{eq:bar_Sigma_k1}
    \begin{split}
      & \bar\Sigma_{12}^{k+1} = (\bar\Sigma_{21}^{k+1})^\top =  
      b_k \,  
      E\left[ \begin{pmatrix} \bar \beta_0 \\
      \bar \gamma_0 \end{pmatrix} 
      \left( M^{k+1}_{B,11} \bar \beta_0 + 
      M^{k+1}_{B,12} \bar \gamma_0 \right) \right] \, \in \reals^2, \\
    & \bar\Sigma_{22}^{k+1} = 
        \delta b_k^2 \,   \left( E[ (M^{k+1}_{B,11} \bar{\beta}_0 + M^{k+1}_{B,12} \bar{\gamma}_0)^2] \, + \, E[g(Z_1^k,   Z_1 + \sigma(Z_2) \bar\varepsilon; \, b_k)^2] \right) \, \in \reals.
\end{split}
\end{equation}

We can now state a convergence result for  $\check \beta^{k+1} \coloneq(\delta b_k) B^{k+1}_1$, which can be viewed as an estimator of $\tilde\beta_0$, the regression parameter in whitened coordinates. For $k \ge 0$, let
\begin{equation}\label{eq: coefficient per iteration}
\mu_{k+1} = (\delta b_k)  M_{B, 11}^{k+1},\qquad \nu_{k+1} = (\delta b_k) M_{B, 12}^{k+1}, \qquad  \tau^2_{k+1}  =  (\delta b_k)^2 T_{B, 11}^{k+1}.
\end{equation}
  The following proposition specifies the asymptotic empirical distribution of $\check \beta^{k+1}$.
\begin{proposition}
    [Finite-iteration law in whitened coordinates]\label{prop: state evolution of beta0}
Suppose that Assumptions~\ref{cond: error distribution}--\ref{cond:white_noise_design}, \ref{cond: empirical distribution of signal matrix}, and \ref{cond: existence initialization} hold. Let $\check \beta^{k+1} = \delta b_k B^{k+1}_1$, where the $B^{k+1}_1$ is the Huber AMP iterate in \eqref{eq: Huber_AMP}. Then, for any continuous function $\phi: \mathbb R^3 \to \mathbb R$ satisfying the polynomial growth condition,
\begin{eqnarray}\label{eq: main theorem  result equation}
\frac{1}{p} \sum_{j =1}^p \phi (\check\beta^{k+1}_j - \mu_{k+1} \tilde\beta_{0, j} - \nu_{k+1} \tilde\gamma_{0, j}, \, \tilde\beta_{0, j}, \,  \tilde\gamma_{0, j} ) \stackrel{}{\to} E [\phi (\tau_{k+1} G, \, \bar\beta_0, \,  \bar\gamma_0)],    
\end{eqnarray}
almost surely as $n, p\to \infty$ and $n/p \to \delta$, where $G \sim N(0,1)$ is independent of $(\bar\beta_0, \bar\gamma_0)$.
\end{proposition}

The proposition is a direct consequence of Lemma~\ref{lem: wasserstein convergence of matrix AMP} and Definition~\ref{def: wasserstein}, which gives an equivalent convergence statement in terms of polynomial growth functions. The coefficients $\mu_{k+1}$ and $\nu_{k+1}$ determine the centring of the whitened-coordinate estimator $\check\beta^{k+1}$, while $\tau_{k+1}^2$ measures the variance from the centre. Taking $\phi(x,b_1,b_2)=x^2$ in Proposition~\ref{prop: state evolution of beta0}  yields, for each $k \ge 0$,
\begin{eqnarray}\label{eq: l2 norm convergence}
\frac{1}{p} \sum_{j =1}^p  (\check\beta^{k+1}_j - \mu_{k+1} \tilde\beta_{0, j} - \nu_{k+1} \tilde\gamma_{0, j})^2 \stackrel{a.s.}{\to} \tau_{k+1}^2 \quad \text{ as } \quad n, p\to \infty.    
\end{eqnarray}

\subsection{Fixed point of Huber AMP}

Let $(b_*, \nu_*, \tau_*^2)$ be an arbitrary solution of the fixed point equations \eqref{eq: fixed point state evolution}--\eqref{eq: barSigma*_def}. We will prove the convergence statement \eqref{eq: asymptotic risk Huber estimator} for this pair $(\nu_*, \tau_*^2)$, and then show that this pair is unique.

We choose the sequence $(b_k)_{k\ge 0}$ such that $b_k  \to b_*$ as $k \to \infty$. The constant $b_* >0$ is specified later.
Then, the fixed point $(B^*_1, \Theta^*_1)$ of the iteration \eqref{eq: Huber_AMP} as $k \to \infty$, if it exists, satisfies:
\begin{equation}
\label{eq: Huber_AMP_FP}
    \begin{split}
        \Theta^*_1 &=  \delta b_* \tbX B^*_1  \, \, - b_* g(\Theta^{*}_1, \bm Y ; \, b_*), \\
        B^{*}_1 & = \tbX^\top g(\Theta^*_1, \bm Y ; \, b_*) \,  -  \,  \delta b_* E[\partial_1 g(Z_1^*, \,  Z_1 + \sigma(Z_2) \bar\varepsilon; \, b_*)] B^*_1.
    \end{split}
\end{equation}
From  \eqref{eq:MBTB_k1}-\eqref{eq:bar_Sigma_k1}, the corresponding state evolution fixed point equations are:
\begin{equation}
\label{eq: state evolution FP MT}
\begin{split}
    & M^*_{B, 11} =   E[ \partial_2 g(Z_1^*,   Z_1 + \sigma(Z_2) \bar\varepsilon; \, b_*) ]   , \quad M^*_{B, 12}=  E[ \partial_2 g(Z_1^*,   Z_1 + \sigma(Z_2) \bar\varepsilon; \, b_*) \sigma'(Z_2) \bar\varepsilon] ,  \\
    & T^*_{B, 11} = E[g(Z_1^*,   Z_1 + \sigma(Z_2) \bar\varepsilon; \, b_*)^2],
\end{split}
\end{equation}
and $(Z_1, Z_2, Z_1^*)^\top \sim N(0, \bar{\Sigma}^*)$, where 
$ \bar \Sigma^* = \begin{pmatrix}
    \bar\Sigma_{11} & \bar\Sigma_{12}^* \\ (\bar{\Sigma}^*_{12})^{\top} &  \bar\Sigma_{22}^*
\end{pmatrix} $, 
with  $\bar\Sigma_{11} = \frac{1}{\delta}E[\bar B \bar B^\top]$, and 
\begin{equation}
\label{eq: Sigma_12_22_FP}
        \begin{split}
      & \bar\Sigma_{12}^{*} = (\bar\Sigma_{21}^{*})^\top =  b_*       E\left[ \begin{pmatrix} \bar \beta_0 \\
      \bar \gamma_0 \end{pmatrix} 
      \left( M^{*}_{B,11} \bar \beta_0 + 
      M^{*}_{B,12} \bar \gamma_0 \right) \right], \\
    & \bar\Sigma_{22}^{*} = 
        (b_*^2 \delta) \left( E[ (M^{*}_{B,11} \bar{\beta}_0 + M^{*}_{B,12} \bar{\gamma}_0)^2] \, + \, E[g(Z_1^*,   Z_1 + \sigma(Z_2) \bar\varepsilon; \, b_*)^2] \right). \\
\end{split}
\end{equation}

We now choose $ b_* = {1}/{(\delta \, M_{B,11}^*)}$, and denote $\nu_* = (\delta b_*) M_{B, 12}^*$ and $  \tau^2_*  = (\delta b_*)^2 T_{B, 11}^*$. We can then rewrite the state evolution fixed point equations \eqref{eq: state evolution FP MT}-\eqref{eq: Sigma_12_22_FP} in terms of $(b_*, \nu_*, \tau^2_*)$, thereby observing that \eqref{eq: state evolution FP MT}-\eqref{eq: Sigma_12_22_FP} are identical to \eqref{eq: fixed point state evolution}-\eqref{eq: barSigma*_def}.  

Next, consider \eqref{eq: Huber_AMP_FP}.  Using  $\delta b _* = {1}/{M_{B,11}^*}$ in the second line and noting from \eqref{eq: g_1_derivatives} that $\partial_1 g = -\partial_2 g$, we obtain that the fixed point of the Huber AMP iteration, if it exists, satisfies:
\begin{equation}
    \label{eq: Huber_AMP_FP_simple}
    \Theta^*_1 = \delta b_* \tbX B^*_1  \, \, - b_* g(\Theta^{*}_1, \bm Y ; \, b_*), \qquad 
    \tbX^\top g(\Theta^{*}_1, \bm Y ; \, b_*) =0.
\end{equation}

The following lemma, proved in Appendix \ref{appendix: proof of AMP FP to Huber}, shows that any fixed point of the Huber AMP corresponds to a solution of \eqref{eq: M-estimator minimization}, i.e., a minimizer of the Huber loss.
\begin{lemma}\label{lemma: AMP FP to Huber}
Let $( B^*_1, \Theta^*_1, b_*)$ be a fixed point of the Huber AMP, i.e., satisfying  \eqref{eq: Huber_AMP_FP_simple}. Let $ \check\beta^* = (\delta b_*) B^*_1$. Then $\hat\beta^* = \Omega_p^{-1/2}\check\beta^*$ is a minimizer of \eqref{eq: M-estimator minimization}.
\end{lemma}

Lemma~\ref{lemma: AMP FP to Huber} motivates the use of the matrix AMP algorithm to analyse the Huber estimator. The result is not used in the subsequent proofs but clarifies the connection between $\hat\beta^*$ and $\hat\beta$.

\subsection{Stationary version of Huber AMP}\label{ssec: stationary Huber AMP}

We now consider a `stationary' version of the Huber AMP iteration in \eqref{eq: Huber_AMP}, defined via the fixed points of the state evolution equations, and a special initialization. Let $(b_*, \nu_*, \tau_*^2)$ be a solution of \eqref{eq: fixed point state evolution}-\eqref{eq: barSigma*_def} (which exists by Proposition \ref{prop: FP existence and uniqueness}). Then, for all $k \ge 0$, we choose
$b_k = b_*$. Using this choice in \eqref{eq: gk_fk_Huber},  for $k \ge 0$ we obtain
\begin{equation} \label{eq: gk_fk_Huber_stationary}
  g_k(w) = g_*(w) \coloneq \begin{pmatrix}
         g(w_1, y ; \, b_*) \\ 0
    \end{pmatrix}, \quad    
    f_{k+1}(w) = f_*(w) \coloneq \begin{pmatrix}
          (\delta b_*) w_1 \\ 0
    \end{pmatrix}, \quad w =(w_1, w_2)^\top \in \reals^2. 
\end{equation}
The AMP iteration \eqref{eq: generic AMP algorithm} is initialized with $\hat{R}^{-1}=0$ and 
\begin{equation}
\label{eq: stationary init}
    \check B^0 =   \begin{pmatrix}
    \tilde \beta_0 \, + \nu_* \tilde \gamma_0 + \, \tau_* \bm v , \ 
    \bm v'
\end{pmatrix}, \qquad \bm v, \, \bm v' \stackrel{\text{i.i.d.}}{\sim} N(0, I_p),
\end{equation}
where the standard Gaussian vectors $\bm v, \bm v'$ are independent of $\tilde \beta_0$. (This is a valid initialization since we only use the AMP iteration as a proof technique, rather than a practical algorithm for estimating $\beta_0$.)  

With the choices in \eqref{eq: gk_fk_Huber_stationary}--\eqref{eq: stationary init}, the state evolution parameters remain stationary at the fixed point values. Indeed, using these choices in \eqref{eq:MBTB_k1}--\eqref{eq:bar_Sigma_k1} together with Assumptions \ref{cond: empirical distribution of signal matrix} and \ref{cond: existence initialization} (stated in Appendix \ref{ssec: assumptions}),  we obtain for $k \ge 1$,
\begin{equation}
    \label{eq: SE stationary}
    \begin{split}
          \bar\Sigma^k = \bar{\Sigma}^*, \quad 
    M^{k+1}_{B,11} = M^{*}_{B,11}, \quad  M^{k+1}_{B,12} = M^{*}_{B,12},   \quad T^{k+1}_{B, 11} = T^*_{B, 11}.
    \end{split}
\end{equation}
The proof of Theorem \ref{thm: dist convergence of huber mean} is now completed via two lemmas. The first, Lemma~\ref{lem: vanishing increments}, shows that for the AMP defined via the stationary functions in \eqref{eq: gk_fk_Huber_stationary}, the iterates satisfy a `Cauchy'-type convergence property, that is, the successive AMP updates become asymptotically indistinguishable after taking the high-dimensional limit and then increasing the iteration number. 

\begin{lemma}[Vanishing increments of stationary AMP recursion]\label{lem: vanishing increments}
Consider the AMP iteration \eqref{eq: generic AMP algorithm}, with the functions and initialization given by \eqref{eq: gk_fk_Huber_stationary}-\eqref{eq: SE stationary}. For $k \ge 1$, let  $\check\beta^k: = \check B^k_1$ and $R_1^k=\bm Y-\Theta_1^k$.  Then, we almost surely have 
\begin{equation}\label{eq: increment vanish beta}
\lim_{k \to \infty} \lim_{n,p \to \infty}
\frac{1}{p}\|\check\beta^{k+1}-\check\beta^k\|_2^2=0, 
\end{equation}
and
\begin{equation}\label{eq: increment vanish residual}
\lim_{k \to \infty}  \lim_{n,p \to \infty}
\frac{1}{n}\|R_1^{k+1}-R_1^k\|_2^2=0,
\end{equation}
where the limits in $n,p$ are taken with $n/p \to \delta$.
\end{lemma} 
The lemma is proved in Appendix \ref{appendix: proof of vanishing increments}. The final lemma shows that $\check\beta^{k+1}$ converges to $\check \beta$ in mean squared error as $k \to \infty$. In addition to Lemma \ref{lem: vanishing increments}, we restrict our attention to Huber objectives $L_n(\beta)$ in \eqref{eq: minimization whitened coordinate} indexed by $a$ satisfying the following restricted secant condition  \citep{polyak1963gradient,zhang2013gradient,guille2022gradient}.

\begin{enumerate}[label=(S6)]
 \item \label{cond: restricted secant inequality}
Consider the gradient of the objective function denoted by $\nabla_\beta L_n (\beta)$. For  the AMP iterate $\check\beta^{k+1}$ and a deterministic constant $\rho>0$, define
\[\kappa(\rho)=\inf_{0<\|h\|_2\leq\rho\sqrt p}
\frac{\left\langle\nabla L_n(\check\beta^{k+1} +h) -\nabla L_n (\check\beta^{k+1}), h \right\rangle}{\|h\|_2^2}.\]
There exist deterministic $\rho >0$ and $\underline{c}>0$
such that
\[
\mathbb P(\liminf_{k\to\infty}\liminf_{p\to\infty}
\kappa(\rho)\geq \underline{c}) = 1.\]
    \end{enumerate}

 To  justify Assumption \ref{cond: restricted secant inequality},  using the formula for $\ell_a$ we write 
    \begin{align}\label{eq: eigenvalue rewrite}
           \left\langle\nabla L_n(\check\beta^{k+1} +h) -\nabla L_n (\check\beta^{k+1}), h \right\rangle= h^\top \left(\sum_{i=1}^n  (\tilde X_i \tilde X_i^\top) \int_0^1 I\left\{|Y_i -  \tilde X_i^\top (\check\beta^{k+1} + th)| \leq a \right\}  \operatorname{d} t \right) h .
    \end{align}
    The above expression suggests that the $i$th residual near the AMP iterate contributes $(\tilde X_i^\top h)^2$ to the gradient inner product of Huber loss \emph{if} it lies in the quadratic region of the Huber loss along the line segment from $\check\beta^{k+1}$ to $\check\beta^{k+1} + h$. 
    We show in 
 Lemma~\ref{lem: inliers near AMP iterate} that for sufficiently small $\eta>0$, the fraction of observations with $|Y_i-\tilde  X_i^\top\check\beta^{k+1}|<a-\eta$ is asymptotically bounded below by a  constant $q_* (\eta)  > \delta^{-1}$. Since $\delta^{-1} = \lim \frac{p}{n} $, the number of these `inlier' observations therefore exceeds $p$, for sufficiently large $p$ and iterations $k$. Thus, under mild assumptions on the whitened covariates, we expect the matrix within parentheses in \eqref{eq: eigenvalue rewrite} to have full rank $p$. Assumption \ref{cond: restricted secant inequality} formalizes this heuristic argument via a lower bound on the minimum eigenvalue of this matrix, providing strictly positive curvature in the neighbourhood of $\check\beta^{k+1}$, for sufficiently large $k$.

    Assumption \ref{cond: restricted secant inequality} is milder than the uniform strong-convexity condition $\inf_x\ell''(x)>0$ used in \citet[Section~4]{donoho2016high} for the analysis of $M$-estimators in the homoscedastic setting. This also differs from the generalized strong convexity condition used to guarantee the existence of the maximum-likelihood estimator, which is used to study the asymptotics of high-dimensional logistic regression \citep{sur2019modern} and other generalized linear models \citep{arpino2026inferring}.

\begin{lemma}[AMP-to-Huber optimizer transfer under restricted secant inequality] \label{lem: mse convergence Huber and AMP at stationarity}
With the same setting as Lemma \ref{lem: vanishing increments} and under Assumption~\ref{cond: restricted secant inequality}, we almost surely have
\begin{equation} \label{eq: mse convergence Huber and AMP at stationarity}
    \lim_{k \to \infty} \limsup_{p \to \infty} \frac{1}{\sqrt p}\|\check\beta - \check\beta^{k+1}\|_2 = 0,
\end{equation}
where $\check \beta$ is the Huber estimator in whitened coordinates defined in \eqref{eq: minimization whitened coordinate}.
\end{lemma}
The proof is given in Appendix \ref{appendix: mse convergence Huber and AMP at stationarity}. We are now ready to complete the proof of Theorem \ref{thm: dist convergence of huber mean}.

For fixed $k$, let \[R_{p,j}^{k+1}
=\check\beta_j^{k+1}
-\tilde\beta_{0,j}
-\nu_{k+1}\tilde\gamma_{0,j}, \  \mbox{ and } \ \hat P_{p,k+1}
=\frac1p\sum_{j=1}^p
\delta_{(
R_{p,j}^{k+1},
\tilde\beta_{0,j},
\tilde\gamma_{0,j})}.\]
That is, $\hat P_{p,k+1}$ is the row-wise joint empirical distribution of $(R_{p}^{k+1}, \tilde\beta_0, \tilde\gamma_0)$. Similarly, let 
\[R_{p,j}
=\check\beta_j
-\tilde\beta_{0,j}
-\nu_*\tilde\gamma_{0,j}, \  \mbox{ and } \  \hat P_{p}
=\frac1p\sum_{j=1}^p
\delta_{(
R_{p,j},
\tilde\beta_{0,j},
\tilde\gamma_{0,j})}.\]
Let $W_2(\hat P_{p,k+1}, \hat P_{p})$ denote the Wasserstein-2 distance (Definition~\ref{def: wasserstein}) between the empirical measures $\hat P_{p,k+1}$  and $ \hat P_{p}$. Since the AMP algorithm provides a coupling between the two empirical measures, we have 
\begin{equation*}
W_2 (\hat P_{p,k+1}, \hat P_{p})^2 \leq \frac{1}{p} \sum_{j=1}^p (\check\beta_{j}^{k+1} - \check\beta_j)^2. 
\end{equation*}
This, together with Lemma \ref{lem: mse convergence Huber and AMP at stationarity}, implies
\begin{equation}
    \label{eq: wasserstein converge empirical measures}
       \lim_{k \to \infty} \lim_{p \to \infty} \, W_2 (\hat P_{p,k+1}, \hat P_{p}) \stackrel{a.s.}{=} 0.
\end{equation}
Furthermore, denoting by $P^*$ the joint law of $(\tau_{*} G, \bar\beta_0, \bar\gamma_0)$, Proposition~\ref{prop: state evolution of beta0} implies that for each  $k\ge 0$
\begin{equation}
\label{eq: SE W2}
      \lim_{p \to \infty} \, W_2 (\hat P_{p,k+1}, P^*) \stackrel{a.s.}{=} 0,
\end{equation}
where we have used the definitions of $\mu_{k+1}, \nu_{k+1}, \tau_{k+1}$ from  \eqref{eq: coefficient per iteration}, and that the state evolution coefficients do not change across iterations (see \eqref{eq: SE stationary}). Therefore, $\nu_{k+1}=\nu_*, \tau_{k+1}= \tau_*$ for all $k \ge 0$.  Using \eqref{eq: wasserstein converge empirical measures}, \eqref{eq: SE W2} and the triangle inequality then yields 
\begin{equation}
\label{eq: W2 final}
    \lim_{p \to \infty} \, W_2 (\hat P_{p}, P^*) \stackrel{a.s.}{=} 0.
\end{equation}
Since Wasserstein-2 convergence implies weak convergence and convergence of second moments \citep[Theorem 6.9]{villani2009optimal}, by uniform integrability we conclude that \eqref{eq: asymptotic risk Huber estimator} holds for any $\phi$ with polynomial growth of order $2$.

The final step in the proof is to argue that $(\nu_*, \tau_*^2)$ is unique.  Let $(b_*, \nu_*, \tau_*^2 )$ be any solution to \eqref{eq: fixed point state evolution}. From \eqref{eq: asymptotic risk Huber estimator}, the following limits hold with probability one:
\begin{align}
\label{eq: lim check beta}
    \lim_{p \to \infty} \frac{1}{p} (\check \beta - \tilde\beta_0)^\top  \tilde\gamma_0 \, = \, \nu_* E[\bar{\gamma}_0^2] \, , \qquad 
    \lim_{p \to \infty} \frac{1}{p} \| \check \beta \|^2 = E[\bar{\beta}_0^2] + (\nu_*)^2 E[\bar{\gamma}_0^2] +  
    2 \nu_* E[\bar{\gamma}_0 \bar{\beta}_0 ] + \tau_*^2.
\end{align}
Since the limits in \eqref{eq: lim check beta} hold with probability one for any solution of \eqref{eq: M-estimator minimization}, we conclude that $\nu_*$ and $\tau_*^2$ are unique.

\section{Discussion}\label{sec: discussion}

This paper develops a precise asymptotic characterization for Huber regression in the proportional regime. The limiting behaviour can be summarized as 
\[
\check\beta_j \approx \tilde\beta_{0,j} + \nu_* \tilde\gamma_{0,j} + \tau_* G, \qquad G \sim N(0,1),
\]
where $( \nu_*, \tau_*^2)$ are determined jointly by the fixed point of a set of equations in \eqref{eq: fixed point state evolution}. The two components have distinct interpretations: $\nu_*$ measures the nuisance coefficient of the variance index parameter $\gamma_0$, and $\tau_* G$ represents fluctuation generated by high-dimensional estimation. The effect $\nu_*$ can be positive, negative, or zero depending on the effect of the variance scaling. Notably, although the heteroscedastic scale is nonlinear, its first-order directional effect on the estimator is linear in $\tilde\gamma_0$. Further, Condition (C1) for $\nu_*=0$ should not be expected without additional conditions on the error distribution and, more broadly, additional conditions on the loss function.

For the linear regression model, the effect of $\tilde\beta_0$ on the estimator is 1 due to a structural simplification in the analysis. This simplification should not be viewed as universal. Indeed, in a generalized linear model, the limiting estimator may then have three components of the form $\mu_*\tilde\beta_0+\nu_*\tilde\gamma_0+\tau_*G$, where $\mu_*$ need not equal 1. Thus, a generalized linear model can exhibit amplification along the mean direction in addition to the effect induced by heteroscedasticity. Analogous rescaling phenomena have been established for homoscedastic high-dimensional logistic regression by \citet{sur2019likelihood}.

As in the homoscedastic setting, the resulting risk is typically nonzero in the proportional asymptotic regime, where $n,p\to\infty$ with $n/p\to\delta\in(1,\infty)$. Hence, the Huber estimator is generally not consistent in normalized $\ell_2$ loss in the proportional regime, and analogous behaviour can be expected for a broader class of unregularized M-estimators. This motivates studying notions of optimality adapted to a nonvanishing limiting risk, including the direct optimization of estimation or out-of-sample prediction error; see, for example, \citet{bellec2025error}.

A related phenomenon, already present in the homoscedastic proportional regime, is the inflation of the effective variance relative to its classical low-dimensional counterpart. This inflation depends on the aspect ratio $\delta$ and, under the relevant regularity conditions, diverges as $\delta\to 1$ in the homoscedastic case. This variance inflation complicates inference and hypothesis testing: as $\delta$ approaches one, the increasing variance reduces precision and makes nontrivial power harder to attain. The exact characterization of the asymptotic variance may allow us to identify the range of $\delta$ over which reliable inference based on a robust asymptotic variance estimator remains possible. 

Finally, matrix approximate message passing serves here as the proof technique for analysing the Huber estimator, but we expect that the analysis can be extended to regularized versions of the estimator. Matrix AMP also suggests a route towards the joint estimation of $\beta_0$ and $\gamma_0$. Existing approaches that account for a variance index parameter commonly model the conditional variance through estimated residuals; see \citet{gu2016high,daye2012high} for sparse high-dimensional regression and \citet{van2010semiparametric} for the classical fixed-$p$ setting. In the proportional regime, however, in-sample residuals retain first-order estimation and leverage effects, so a naive residual-then-scale procedure need not behave as it does in fixed dimensions. This motivates a joint location-scale objective or a coupled system of robust estimating equations that estimates the mean and variance indices simultaneously. Joint quasi-likelihood estimation offers one possible formulation; see \citet{Antoniadia2016joint,gijbels2019robust}. Extending matrix AMP to analyse such procedures would clarify how estimating the variance index parameter direction affects the risk of the mean estimator and inform the joint choice of any location and dispersion tuning parameters.

Taken together, the results show that heteroscedasticity affects high-dimensional robust regression through both asymptotic variance and a systematic directional displacement. Notions of optimality in this regime must therefore account jointly for signal rescaling, the variance-index effect, and high-dimensional fluctuation.

\subsection*{AI Declaration} The authors used OpenAI ChatGPT 5.5 and 6, and Claude Opus 5.5 to assist with the following: writing code for simulations, the proof of Proposition \ref{prop: FP existence and uniqueness} (to prove boundedness of the function $f_2$),  exploring possible formulations of Assumption~\ref{cond: restricted secant inequality}, and proofreading. 
The authors wrote the paper themselves and take full responsibility for the content. 


\bibliographystyle{abbrvnat_nourl}
\bibliography{heteroscedasticity}

\clearpage
\appendix
\numberwithin{equation}{section}
\numberwithin{definition}{section}
\numberwithin{lemma}{section}
\numberwithin{remark}{section}

\section{Regularity conditions and Auxiliary results}
Assumptions~\ref{cond: scale function}--\ref{cond: empirical distribution of signal matrix} for Theorem \ref{thm: dist convergence of huber mean} have been stated in Section~\ref{sec: Huber mean main section}. We state  an additional Assumption~\ref{cond: existence initialization} here, which is used for Lemma \ref{lem: wasserstein convergence of matrix AMP},  a Wasserstein convergence result for generic matrix AMP. We note that \ref{cond: existence initialization} is not an additional assumption for Theorem \ref{thm: dist convergence of huber mean}.
\subsection{Regularity conditions}\label{ssec: assumptions}
\begin{enumerate}[label=(S7)]
        \item \label{cond: existence initialization} 
    The initialization $\check B^0$ for the recursion \eqref{eq: generic AMP algorithm} satisfies the following. The joint empirical distribution of $(\tilde B, \check B^0)$ converges to the distribution of $(\bar B^\top, (\bar{B}^0)^\top)^\top$ in Wasserstein-2 distance  as $n, p\to \infty$ with $n/p \to \delta > 0$. The random vector $(\bar B^\top, (\bar{B}^0)^\top)^\top \in \reals^4$ is assumed to have a positive definite covariance matrix $\bar\Sigma^0$. This implies that 
    \begin{eqnarray}\label{eq: initialization sigma exist}
       \frac{1}{n} \begin{pmatrix}
       \tilde B^\top \tilde B & \tilde B^\top  \check B^0\\
       (\check B^0)^\top \tilde B & (\check B^0)^\top (\check B^0)
       \end{pmatrix} \stackrel{a.s.}{\to}  \bar\Sigma^0 = 
       \begin{pmatrix}
       \bar\Sigma_{11} & \bar\Sigma_{12}^0 \\
    \bar\Sigma_{21}^0 & \bar\Sigma_{22}^0
       \end{pmatrix}, 
    \end{eqnarray}
    where $\bar \Sigma^0 \in \mathbb R^{4 \times 4}$ is  positive definite  with  $\bar\Sigma_{11} = \frac{1}{\delta} E[\bar{B} \bar{B}^\top]$ (the signal strength for $\bar B$ in Assumption~\ref{cond: empirical distribution of signal matrix}). 
\end{enumerate}

\subsection{Auxiliary results} \label{ssec: aux_results}

\begin{definition}[Wasserstein distance \citep{panaretos2019statistical}]\label{def: wasserstein distance}
    Let $\mathcal P_r(\mathbb R^k)$ denote the set of probability measures on $\mathbb R^k$ with finite $r$th moments for $r \in [1,\infty)$. For $P,Q\in\mathcal P_r(\mathbb R^k)$, define
\[W_r(P,Q)=
\left\{
\inf_{\pi\in\Pi(P,Q)}
\int
\|x-y\|_2^r\,\pi(dx,dy)
\right\}^{1/r},
\]
where $\Pi(P,Q)$ is the set of couplings of $P$ and $Q$.
\end{definition}

\begin{definition}[Wasserstein convergence \citep{wang2024universality}]\label{def: wasserstein}
For $n \in \mathbb N$ and $r \in [1, \infty)$, define vectors $a^1, \ldots, a^k \in \mathbb R^n$ and a random vector $(A^1, \ldots, A^k) \in \mathbb R^k$, we write \[
(a^1, \ldots, a^k) \stackrel{W_r}{\to} (A^1, \ldots, A^k), n\to \infty,
\]
for the Wasserstein-$r$ convergence of the empirical distribution of the rows of $(a^1, \ldots, a^k) \in \mathbb R^{n\times k}$ to the joint law of $(A^1, \ldots, A^k)$. This means, for any continuous function $\phi: \mathbb R^k \to \mathbb R$ and input vector $(a_i^1, \ldots, a_i^k) \in \mathbb R^k$ satisfying the polynomial growth condition (of order $r$)
\begin{eqnarray}\label{eq: polynomial growth}
    |\phi (a_i^1, \ldots, a_i^k)| \leq C (1 + \|(a_i^1, \ldots, a_i^k )\|_2^r)  \mbox{ for a constant } C > 0,
\end{eqnarray}
we have as $n \to \infty$,
\[
\frac{1}{n} \sum_{i =1}^n \phi(a_i^1, \ldots, a_i^k) \to E [\phi (A^1, \ldots, A^k)] .
\]
\end{definition}

The following convergence result for the iterates of the matrix AMP in \eqref{eq: generic AMP algorithm} follows directly from \cite[Theorem 3.4]{tan2024approximate}.
\begin{lemma}[Wasserstein convergence of generic matrix AMP]\label{lem: wasserstein convergence of matrix AMP}
Consider the recursion in \eqref{eq: generic AMP algorithm} for the model in \eqref{eq: data generating process}, with the whitened matrix $\tilde{\bm X}= \bm X \Omega_p^{-1/2}$ satisfying Assumption~\ref{cond:white_noise_design}. Further, assume that Assumptions~\ref{cond: empirical distribution of signal matrix}, \ref{cond: error distribution}, and \ref{cond: existence initialization} hold. Then for each $k \ge 0$,
\begin{align}
    &(\tilde B, \check B^0, \, B^1, \ldots, B^{k+1}) \stackrel{W_2}{\to} 
(\bar B, \bar B^0  , M_B^{1}\bar B + G _B^{1}, \ldots,  M_B^{k+1}\bar B + G _B^{k+1})  \label{eq: wasserstein convergence 1}, \\
 & (\bm\varepsilon, \tilde\Theta, \Theta^0, \Theta^1, \ldots, \Theta^{k}) \stackrel{W_2}{\to} (\bar\varepsilon, Z, Z^0, Z^1, \ldots, Z^{k})  \label{eq: wasserstein convergence 2},
\end{align}
almost surely as $n, p \to \infty$ and $n/p \to \delta$. In the above, for $t=0,\ldots, k$, $G_B^{t+1} \sim N(0, T_B^{t+1})$ is independent of $\bar B$, and  $\begin{pmatrix} Z \\ Z^t \end{pmatrix} \sim N(0, {\bar\Sigma}^t)$, where the matrices 
$M_B^{t+1}, T_B^{t+1}, {\bar\Sigma}^t$ are given by the state evolution recursion in \eqref{eq: state-evolution matrix AMP}-\eqref{eq:Sigmak_update}, with initialization $\bar\Sigma^0$.
Further, for $0 \le s \le s' \le k$, we have
\begin{align}
& E [G_B^{s+1} (G_B^{s'+1})^\top]   = E[h_s(Z, Z^s, \bar\varepsilon) h_{s'}(Z, Z^{s'}, \bar\varepsilon)^\top], \label{eq: GB_covariance} \\
& E [Z^s (Z^{s'})^\top]  = \frac{1}{\delta} E[ f_{s}(M_B^{s} \bar{B} + G_B^{s}) (f_{s'}(M_B^{s'} \bar{B} + G_B^{s'}))^\top].  \label{eq: Gtheta_covariance}
\end{align}
(In \eqref{eq: Gtheta_covariance}, for $s=0$, we take $f_{0}(M_B^{0} \bar{B} + G_B^{0}) \equiv \bar{B}^0$.)
\end{lemma}

\begin{remark} \label{rem: det_vs_empirical}
    In the matrix AMP analysed in \cite{tan2024approximate}, the `Onsager' matrices $\bar{F}^k$ and $\bar{C}^k$ in \eqref{eq: CF limiting} are defined via empirical averages rather than expectations. Lemma \ref{lem: wasserstein convergence of matrix AMP} holds for the matrix AMP iteration \eqref{eq: generic AMP algorithm} with either empirical or deterministic versions of $\bar{F}^k, \bar{C}^k$. In fact, the proof of \cite[Theorem 3.4]{tan2024approximate} is based on a reduction to an abstract AMP iteration analysed by \cite{wang2024universality}, which uses deterministic Onsager coefficients.
\end{remark}

 \begin{lemma}[Inliers near the stationary AMP iterate]\label{lem: inliers near AMP iterate}  Consider the setting of Lemma \ref{lem: vanishing increments}.
        For $a \in (0, \infty)$ and $\eta \in (0, a/2)$, let $q_* (\eta) =  \mathbb P\left( |Z_1 + \sigma(Z_2)\bar\varepsilon - Z_1^* | < (a - 2\eta)(1 + b_*) \right)$, then it holds that $\lim_{\eta \to 0}  q_*(\eta) =  \delta^{-1} (1 + b_*^{-1})$. Further, define
        $q_k (\eta) = \frac{1}{n} \sum_{i=1}^n I\left\{ |Y_i - \tilde X_i^\top \check\beta^{k+1}| < (a - \eta) \right\}$, it holds that 
        \[\mathbb P (\liminf_{k\to\infty}\liminf_{n, p \to\infty} q_k (\eta ) \ge q_*(\eta) ) = 1 . \]
    \end{lemma}

The proof of Lemma~\ref{lem: inliers near AMP iterate} is given in Appendix~\ref{ssec: appendix proof of sufficient inliers}.

\section{Proofs }

\subsection{Proof of Proposition \ref{prop: FP existence and uniqueness} (Existence  of solution to fixed point equations)} \label{appendix: proof FP existence and uniqueness}
\begin{proof}
For brevity, we drop the subscripts $*$ throughout the proof. We will assume  that $a <\infty$, as the case of $a=\infty$ (least squares) is covered by Example \ref{example: LS}.
We first show that for any  $(\nu, \tau^2)$, the solution for $b$ in the first equation of \eqref{eq: fixed point state evolution} exists and is unique.  Using \eqref{eq: g_1_derivatives}, this equation can be written as
\begin{align}
\label{eq: FP equation b}
    \frac{1}{b} +1  = \delta \,  \mathbb P( |Z_1 + \sigma(Z_2) \bar\varepsilon - Z_1^*|\le a(1+b)) , 
\end{align}
 where we note from \eqref{eq: barSigma*_def} that the joint law of $(Z_1, Z_2, Z_1^*)$ depends only on $(\nu, \tau^2)$, and not on $b$. For $b \in (0, \infty)$, the left hand side of \eqref{eq: FP equation b} strictly decreases from $\infty$ to $1$, while the right hand side increases from $\delta P( |Z_1 + \sigma(Z_2) \bar\varepsilon - Z_1^*|\le a) \ge 0$ to $\delta$. Therefore, \eqref{eq: FP equation b} has a  unique solution $b(\nu, \tau^2)$  if (and only if) $\delta >1$. Moreover, for any $\tau^2>0$, the right side of \eqref{eq: FP equation b} is continuous in $b$, and therefore the solution $b(\nu, \tau^2)$ is continuous in $(\nu, \tau^2)$.

Writing $b \equiv b(\nu, \tau^2)$, we can express the  set of fixed point equations in \eqref{eq: fixed point state evolution} in terms of  $(\nu, \tau^2)$. Let
$f \equiv (f_1, f_2) :  \reals \times \reals_+ \to \reals \times \reals_+$, where
\begin{align}
    f_1(\nu, \tau^2) & =  \frac{ \delta b }{1+b} \,  E\big[ \, I \big\{ |Z_1 + \sigma(Z_2) \bar\varepsilon - Z_1^*| \le a \big(1+ b \big)  \big\} \, \sigma'(Z_2)\bar\varepsilon \big]  \, ,  \label{eq:f1_def} \\
    f_2(\nu, \tau^2) & =  (\delta b)^2 E\Big[g (Z_1^*,  Z_1 + \sigma(Z_2) \bar\varepsilon ; \, b)^2\Big] .
    \label{eq:f2_def}
\end{align}
We first argue that $f \equiv (f_1, f_2)$ is continuous in both $\nu$ and  $\tau^2$. Indeed, the expectations in \eqref{eq:f1_def}, \eqref{eq:f2_def} are both bounded (by Assumption
\ref{cond: scale function} and \eqref{eq: g_1 function Huber}, respectively), and continuity follows by dominated convergence and the continuity of $b$. 

Next we show that there exists a compact convex set $\mathcal K$ of the form $$\mathcal K \equiv  [-K , \, K] \times [\underline{A}, \, A], $$ such that $f(\nu, \tau^2) \in \mathcal K$ for all $(\nu, \tau^2) \in \mathcal K$.  Recalling that $f$ is continuous on $\mathcal K$, the existence of a fixed point is then guaranteed by Brouwer's fixed point theorem \citep{zeidler1988nonlinear}.

The constant $K >0$  can be taken to be $\delta E[\bar\varepsilon^2]^{1/2} (E[|\sigma'(Z_2) |^2])^{1/2}$, which is finite by Assumptions~\ref{cond: scale function} and \ref{cond: error distribution}. Then, for $f_1$, it holds by the Cauchy-Schwarz inequality that $|f_1(\nu, \tau^2)| < \delta E[|\sigma'(Z_2)\bar\varepsilon|] 
\leq \delta E[\bar\varepsilon^2]^{1/2}{E[|\sigma'(Z_2)|^2]}^{1/2} = K$. The remainder of the proof consists of showing the existence of constants  $A> \underline{A} >0$, completing the specification of $\mathcal K$. Let 
\begin{align*}
    W= (Z_1 + \sigma(Z_2) \bar\varepsilon - Z_1^*), \quad  \text{and} \quad T_{a(1+b)}(W) = \begin{cases} W, & |W|  \le a(b+1),  \\
    a(b+1) \cdot \text{sgn}(W), & \text{ otherwise}.\end{cases}
\end{align*}
Then, using \eqref{eq: g_1 function Huber} we can  write $f_2$ as 
\begin{align}
\label{eq:f2_T_def}
    f_2(\nu, \tau^2) = \frac{(\delta b)^2}{(1+b)^2} E\left[ T_{a(1+b)}(W)^2 \right]
    =  \frac{E\left[ T_{a(1+b)}(W)^2 \right]}{\mathbb P( |W| \le a(1+b))^2} \, ,
\end{align}
where the last equality uses \eqref{eq: FP equation b}.  Recalling that $(Z_1, Z_2, Z_1^*) \sim N(0, \bar\Sigma^*)$, using the covariance in \eqref{eq: barSigma*_def} we have $$W \stackrel{d}{=} - \nu Z_2 \,  +  \sigma(Z_2) \bar\varepsilon \, +  \frac{\tau}{\sqrt{\delta}} G, \qquad Z_2 \perp G \sim N(0,1).  $$
Denoting $U=- \nu Z_2 \,  +  \sigma(Z_2) \bar\varepsilon$, we claim that 
\begin{align}
\label{eq:Tvar_decomp}
   E\left[ T_{a(1+b)}(W)^2 \right] \le \frac{\tau^2}{\delta} \mathbb P( |W| \le a(1+b)) + E[U^2].
\end{align}
To show \eqref{eq:Tvar_decomp},  consider the conditional second moment
\begin{align}
\label{eq:cond_second_mom_decomp}
    E\left[ T_{a(1+b)}(W)^2 \mid U \right] = \operatorname{var}\Big( T_{a(1+b)}(U + \tau/\sqrt{\delta} \, G) \mid U \Big) + E[T_{a(1+b)}(U + \tau/\sqrt{\delta} \,  G ) \mid U]^2. 
\end{align}
The function $g \mapsto T_{a(1+b)}(U + \tau/\sqrt{\delta} \, g)$ is Lipschitz with (weak) derivative $\frac{\tau}{\sqrt{\delta}}I(|U + \tau/\sqrt{\delta} \, g| \le a(b+1))$. Therefore, by the Gaussian Poincar{\'e} inequality \citep[Theorem 3.20]{boucheron2013concentration},
\begin{align}
\label{eq:var_poincare_bound}
    \operatorname{var}\Big( T_{a(1+b)}(U + \tau/\sqrt{\delta} \, G) \mid U \Big) \le \frac{\tau^2}{\delta} \mathbb P(|U + \tau/\sqrt{\delta} \, G| \le a(1+b)).
\end{align} 
Since $T_{a(1+b)}$ is 1-Lipschitz, so is $\Psi(U) \coloneq  E[T_{a(1+b)}(U + \tau/\sqrt{\delta} \,  G ) \mid U]$, which implies $|\Psi(U) -\Psi(0)| \le |U|$.  Moreover, since $T_{a(1+b)}$ is odd, $\Psi(0)=0$, hence 
\begin{align}
\label{eq:second_mom_T_bound}
  \Psi(U)^2 =  E[T_{a(1+b)}(U + \tau/\sqrt{\delta} \,  G ) \mid U]^2  \le U^2. 
\end{align}
Using  \eqref{eq:var_poincare_bound} and \eqref{eq:second_mom_T_bound} in   \eqref{eq:cond_second_mom_decomp} and taking expectations yields \eqref{eq:Tvar_decomp}. Using \eqref{eq:Tvar_decomp} in \eqref{eq:f2_T_def}, we obtain
\begin{align}
    f_2(\nu, \tau^2) & \le \frac{\tau^2}{\delta \, \mathbb P(|U + \tau/\sqrt{\delta} \, G| \le a(1+b))} \, + \, \frac{E[U^2]}{\mathbb P(|U + \tau/\sqrt{\delta} \, G| \le a(1+b))^2} \nonumber \\
    & =  \tau^2\frac{b}{1+b} \, + \, \frac{E[U^2]}{\mathbb P(|U + \tau/\sqrt{\delta} \, G| \le a(1+b))^2}  \nonumber \\
    & \le \tau^2\frac{b}{1+b} \, + \, \delta^2 E[U^2] = \tau^2 -\frac{1}{1+b}\tau^2 \, + \, \delta^2 E[U^2], \label{eq:f2_bound}
\end{align}
where the first equality follows from \eqref{eq: FP equation b}. Next, we show that $(1+b)$ grows at most linearly in $\tau$, and use it in  \eqref{eq:f2_bound} to obtain $0 <\underline{A} <A$ such that $f_2(\nu, \tau^2) \in [\underline{A}, A]$ for $\tau^2 \in [\underline{A}, A]$ and  $\nu \in [-K, K]$.
Specifically, we show
\begin{align}
\label{eq:b1_bound}
    (1+b) \le \max\left\{ 1+ \frac{2}{\delta-1}, \ \, \frac{1}{a} \sqrt{\frac{2 \delta}{(\delta -1)} \left( \frac{\tau^2}{\delta} + E[U^2]\right)}\right\},
\end{align}
noting that the right hand side is linear in $\tau$. The bound is trivial for $b \le \frac{2}{\delta -1}$. For $b > \frac{2}{\delta -1}$,  from \eqref{eq: FP equation b} we have
\begin{align}
\label{eq:LB_tail_prob}
    P(|W| > a(1+b) ) > \frac{\delta-1}{2 \delta},
\end{align}
and using Chebyshev's inequality,  we have the upper bound
\begin{align}
\label{eq:UB_tail_prob}
   \mathbb P(|W| > a(1+b)) \le \frac{E[W^2]}{a^2(1+b)^2} = \frac{E[U^2] + \frac{\tau^2}{\delta}}{a^2(1+b)^2}.
\end{align}
Combining the bounds in \eqref{eq:LB_tail_prob} and \eqref{eq:UB_tail_prob} yields \eqref{eq:b1_bound}. By \eqref{eq:b1_bound},
\[
  \frac{\tau^2}{1+b} \;\ge\;  \min\Big\{ \tau^2 \frac{(\delta -1)}{(\delta+1)},\; \frac{\tau^2}{\Lambda_2} \Big\}, \qquad \Lambda_2: = \frac{1}{a} \sqrt{\frac{2 \delta}{(\delta -1)} \left( \frac{\tau^2}{\delta} + E[U^2]\right)}.
\]
By the definition of $U$, $E[U^2] = \nu^2 E[Z_2^2] + E[\sigma(Z_2)^2]E[\bar \varepsilon^2] $.
With $K$ chosen as described above, let $\overline{M}^2 \coloneq K^2 E[Z_2^2] + E[\sigma(Z_2)^2]E[\bar\varepsilon^2]$, so that $E[U^2] \le \overline{M}^2$ uniformly over $|\nu|\le K$.  Therefore, for
$|\nu|\le K$ and $\tau^2 \ge  \delta \overline{M}^2$, we have $\Lambda_2  \le \frac{2 \tau}{a \sqrt{\delta -1}}$ and $\tau^2/\Lambda_2\ge \frac{a}{2}\sqrt{\delta-1}\;\tau$. Using this in \eqref{eq:f2_bound}, $f_2(\nu, \tau^2) \le \tau^2$ whenever $| \nu | \le K$ and 
\begin{equation}\label{eq: T choice}
  \tau \;\ge\; \tau_0 \;\coloneq\; \max\left\{ \sqrt{\delta}\,\overline{M}, \;\;
  \delta\,\overline M\sqrt{\frac{\delta+1}{\delta-1}}, \;\;
  \frac{2\delta^2 \overline{M}^2}{a\sqrt{\delta-1}} \right\} .
\end{equation}

We set $A\coloneq \tau_0^2 + \delta^2 \overline{M}^2$, and note from \eqref{eq:f2_bound} that $f_2(\nu, \tau^2) \le A$ for $\tau^2 \le \tau_0^2$. For $\tau^2 \in [\tau_0^2 , A]$, by \eqref{eq: T choice} we have $f_2(\nu, \tau^2) \le \tau^2 \le A$ . To obtain a lower bound $\underline{A}$,  we start from \eqref{eq:f2_T_def}, which implies $f_2(\nu, \tau^2) \ge E[T_{a(1+b)}(W)^2]$. Note that  the map  $b \mapsto T_{a(1+b)}(W)^2$ is non-decreasing in $b$ for each $W$, and from \eqref{eq: FP equation b},  $b \ge \frac{1}{\delta -1}$. Therefore,
\[ f_2(\nu, \tau^2) \ge E[T_{a(1+b)}(W)^2] \ge E[T_{a\delta/(\delta-1)}(W)^2]. \] 
The right side is
bounded below by a constant $\underline{A} > 0$ depending only on $a,\delta$ and the laws of $Z_2,\bar\varepsilon$. Shrinking $\underline{A}$ if necessary so that $A >\underline{A}$ completes the specification of the compact convex set $\mathcal K \equiv  [-K , \, K] \times [\underline{A}, \, A]$ such that $f(\nu, \tau^2) \in \mathcal K$ for all $(\nu, \tau^2) \in \mathcal K$. This completes the proof.
 \end{proof}

\subsection{Proof of Lemma~\ref{lemma: AMP FP to Huber}}\label{appendix: proof of AMP FP to Huber}

\begin{proof}
By the definition of $\bar g$ in \eqref{eq: proximal}, for any $(u, y) \in \reals^2$ and $b \in \reals_+$, we have
\begin{eqnarray}\label{eq: proximal estimating equations}
    \partial_1 \ell_a (\bar g (u, y; \, b), y) + \frac{1}{b} (\bar g(u, y; \, b) - u) =0.
\end{eqnarray}
At the fixed point,  by \eqref{eq: Huber_AMP_FP_simple} we have
\begin{eqnarray}\label{eq: X orthogonal to residual}
 {\bm X}^\top g(\Theta^{*}_1, \bm Y ; \, b_*) = (\tilde{\bm X} \Omega_p^{1/2})^\top g(\Theta^{*}_1, \bm Y ; \, b_*) = \Omega_p^{1/2} \tilde{\bm X}^\top g(\Theta^{*}_1, \bm Y ; \, b_*) =  0.   
\end{eqnarray}
Further, by  \eqref{eq: Huber_AMP_FP_simple} and \eqref{eq: effective score}, 
\begin{eqnarray}\label{eq: proximal at stationarity}
  \bar g(\Theta_1^*, \bm Y; \, b_*) = \tilde{\bm X}\check\beta^* = \bm X \hat\beta^* ,
\end{eqnarray}
where we have used the definitions of $\check\beta^* $ and $\hat\beta^*$ in the lemma statement. Therefore,
\begin{align*}
   \sum_{i=1}^n \partial_1 \ell_a(X_i^\top\hat\beta^*, Y_i)X_i &=  \sum_{i=1}^n \partial_1 \ell_a( \bar g(\theta_{1i}^*, Y_i; \, b_*), Y_i) X_i \\
    & = - \sum_{i=1}^n g(\theta_{1i}^*, Y_i; \, b_*) X_i = - \bm X^\top g(\Theta^{*}_1, \bm Y ; \, b_*)  = 0_p,
\end{align*}
where the first equality holds by \eqref{eq: proximal at stationarity}, the second  by \eqref{eq: proximal estimating equations}, and the last by \eqref{eq: X orthogonal to residual}. 

Thus, it holds that
\[
\nabla_\beta \sum_{i=1}^n \ell_a( X_i^\top \beta, Y_i) \Big\vert_{\, \beta = \hat\beta^*} = \sum_{i=1}^n \partial_1 \ell_a(X_i^\top \beta, Y_i)X_i \Big\vert_{\, \beta = \hat\beta^*}  =  0_p,
\]
so convexity of $\ell_a$ implies that $\hat\beta^*$ is a minimizer of \eqref{eq: M-estimator minimization}.
\end{proof}

\subsection{Proof of Lemma \ref{lem: vanishing increments}} \label{appendix: proof of vanishing increments}
\begin{proof}
We write
\begin{align}
\label{eq: beta three term split}
    \frac{1}{p}\|\check\beta^{k+1}-\check\beta^k\|_2^2 =
\frac{1}{p}\|\check\beta^{k}\|_2^2 \,  + \,  \frac{1}{p}\|\check\beta^{k+1}\|_2^2 
- 2 \frac{1}{p}(\check\beta^{k})^\top \check\beta^{k+1},  \\
\label{eq: R three term split}
\frac{1}{n}\|R_1^{k+1}-R_1^k\|_2^2 =   
\frac{1}{n}\|\Theta_1^k\|_2^2 \,  + \,  \frac{1}{n}\|\Theta_1^{k+1}\|_2^2 
- 2 \frac{1}{n}(\Theta_1^{k})^\top \Theta_1^{k+1}.
\end{align}
We compute the limiting values for each of the terms on the right using the state evolution result in Lemma \ref{lem: wasserstein convergence of matrix AMP} (and the special case in Proposition \ref{prop: state evolution of beta0}). The limits in the rest of the proof all hold with probability $1$.

For the first two terms in \eqref{eq: beta three term split}, we apply Proposition~\ref{prop: state evolution of beta0} with the test function $\phi: \reals^3 \to \reals$, defined as $\phi(u,v,w) = (u + \mu_{k+1}  v+ \nu_{k+1} w)^2$. This yields
\begin{equation}
\lim_{p \to \infty} \frac{1}{p}\|\check\beta^{k}\|_2^2 = E \left[ (\tau_k G + \mu_{k} \bar \beta_0 + \nu_{k} \bar \gamma_0)^2 \right], \qquad  \lim_{p \to \infty} \frac{1}{p}\|\check\beta^{k+1}\|_2^2 = E \left[ (\tau_{k+1} G + \mu_{k+1} \bar \beta_0 + \nu_{k+1} \bar \gamma_0)^2 \right],
\end{equation}
where $G \sim N(0,1)$ is independent of $(\bar \beta_0, \bar \gamma_0)$.
Using \eqref{eq: coefficient per iteration} and  the stationary functions and state evolution parameters in \eqref{eq: gk_fk_Huber_stationary}-\eqref{eq: SE stationary}, for $k \ge 0$ we have 
\[ \mu_{k} = (\delta b_*) M^*_{B, 11} =1, \qquad \nu_k = (\delta b_*) M^*_{B, 12} = \nu_*, \qquad \tau_k^2 = (\delta b_*)^2 T^*_{B, 11} =\tau^2_*\, . \]
Therefore, 
\begin{equation}
\label{eq: betak norm}
    \lim_{p \to \infty} \frac{1}{p}\|\check\beta^{k}\|_2^2 = E[ (\bar \beta_0 + \nu_{*} \bar \gamma_0)^2] + \tau_*^2, \qquad k \ge 0.
\end{equation}
Similarly, using Lemma \ref{lem: wasserstein convergence of matrix AMP}, we have 
\begin{equation}
\label{eq: thetak norm}
     \lim_{p \to \infty} \, \frac{1}{n}\|\Theta_1^k\|_2^2 = \bar \Sigma^k_{22} = \bar \Sigma^*_{22}, \qquad k \ge 0,
\end{equation}
where the last equality holds due to the stationary state evolution update \eqref{eq: SE stationary}.

We now address the cross terms in \eqref{eq: beta three term split} and \eqref{eq: R three term split}.
Next, recalling that $\check\beta^{k} = (\delta b_*) B^k_1$ and using  \eqref{eq: wasserstein convergence 1} (and Definition \ref{def: wasserstein}), we have
\begin{equation}
\label{eq: Bk_Bk1}
\begin{split}
        \frac{1}{p}(\check\beta^{k})^\top \check\beta^{k+1} & = (\delta b_*)^2 \frac{1}{p} \sum_{j=1}^p  B^k_{1,j}  B^{k+1}_{1,j} \\
        & \stackrel{p \to \infty}{\longrightarrow} \, (\delta b_*)^2 \left(  E[ (M^k_{B, 11} \bar \beta_0 + M^k_{B, 12} \bar \gamma_0 ) (M^{k+1}_{B, 11} \bar \beta_0 + M^{k+1}_{B, 12} \bar \gamma_0 )] \,  +  \, E[G^k_{B,1} G^{k+1}_{B,1}] \right) \\
        & = E[ (\bar \beta_0 + \nu_{*} \bar \gamma_0)^2] \, + \, 
        (\delta b_*)^2 E[G^k_{B,1} G^{k+1}_{B,1}],
\end{split}
\end{equation}
where the last line again holds due to \eqref{eq: SE stationary}. For the cross-term in \eqref{eq: R three term split}, using \eqref{eq: Gtheta_covariance} we have 
\begin{equation}
\label{eq: Zk_Zk1}
\begin{split}
    \lim_{n \to \infty} \frac{1}{n}(\Theta_1^{k})^\top \Theta_1^{k+1} & =    E[Z^k_1 Z^{k+1}_1 ] \\
    & = \frac{1}{\delta }(\delta b_*)^2   E[ (M^k_{B, 11} \bar \beta_0 + M^k_{B, 12} \bar \gamma_0  + G^k_{B,1}) (M^{k+1}_{B, 11} \bar \beta_0 + M^{k+1}_{B, 12} \bar \gamma_0 + G^{k+1}_{B,1} )] \\
    & = \frac{1}{\delta} \left( E[ (\bar \beta_0 + \nu_{*} \bar \gamma_0)^2]  \, + \, 
        (\delta b_*)^2 E[G^k_{B,1} G^{k+1}_{B,1}] \right).
    \end{split}
\end{equation}
From \eqref{eq: GB_covariance} and recalling that $h_k(s, w, v) =  g_k (w, q(s, v))$, we have 
\begin{equation}
\label{eq: GB_k_k1}
\begin{split}
        E[G^k_{B,1} G^{k+1}_{B,1}] & = E[ h_{k-1, 1}(Z, Z^{k-1} , \bar \varepsilon)
        h_{k, 1}(Z, Z^{k} , \bar \varepsilon)] \\
        & = E[g(Z_1^{k-1},   Z_1 + \sigma(Z_2) \bar\varepsilon; \, b_*) g(Z_1^{k},   Z_1 + \sigma(Z_2) \bar\varepsilon; \, b_*)],
\end{split}
\end{equation}
where we used the definition of $g_k$ in \eqref{eq: gk_fk_Huber_stationary}.  

The equations \eqref{eq: Zk_Zk1} and \eqref{eq: GB_k_k1} together define a recursion in the one-step covariance $E[G^k_{B,1} G^{k+1}_{B,1}]$, for $k \ge 0$.
From  \eqref{eq: SE stationary}, we know that marginally $G^k_{B,1} \sim N(0, T^*_{B,11})$ for  all $k \ge 0$, where  $T^*_{B,11}=  \tau_*^2/(\delta b_*)^2$.  Using the initialization \eqref{eq: stationary init} together with Assumption \ref{cond: existence initialization}, we have $E[G^0_{B,1} G^{k}_{B,1}] = 0$ for $k \ge 1$. Defining the correlation coefficient 
$$q_k \coloneq \frac{E[G^k_{B,1} G^{k+1}_{B,1}]}{T^*_{B,11}}  =  \frac{E[G^k_{B,1} G^{k+1}_{B,1}] (\delta b_*)^2}{\tau_*^2} ,$$
we now show that $\lim_{k \to \infty} q_k =1$. Using this in \eqref{eq: Bk_Bk1} and \eqref{eq: Zk_Zk1} and combining with the limits in \eqref{eq: betak norm}-\eqref{eq: thetak norm} yields the lemma.

To show $q_k \to 1$, we begin by writing the  one-step recursion from \eqref{eq: Zk_Zk1} and \eqref{eq: GB_k_k1} as $q_k =H(q_{k-1})$, for $k \ge 1$, where 
\begin{equation}
\label{eq: Hdef}
\begin{split}
    H(q) = 
    (T^*_{B,11})^{-1}E[g(\hat Z_1,   Z_1 + \sigma(Z_2) \bar\varepsilon; \, b_*) g(\check Z_1,   Z_1 + \sigma(Z_2) \bar\varepsilon; \, b_*)], \quad q \in [0,1].
\end{split}
\end{equation}
Here, the expectation is with respect to the following distribution:   $(Z_1, Z_2, \hat Z_1, \check Z_1 )^\top$ are jointly Gaussian with
 $(Z_1, Z_2, \hat Z_1)^\top \stackrel{\rm d}{=} (Z_1, Z_2, \check Z_1)^\top  \sim N(0, \bar{\Sigma}^*)$, and 
 $ E[ \hat Z_1 \check Z_1 ] = \frac{1}{\delta} \left( E[ (\bar \beta_0 + \nu_{*} \bar \gamma_0)^2]  \, + \, 
        q \tau_*^2  \right)$. 
The recursion is initialized with $q_0=0$.  Using the definition of $\bar \Sigma^*$ from \eqref{eq: barSigma*_def} and properties of jointly Gaussian random vectors, we can express $\hat Z_1, \check Z_1$ as:
\begin{align}
\label{eq: hat check exp}
    \hat Z_1 = \bar{M} Z + \,  \bar{\zeta} G_1, \qquad \check Z_1 = \bar{M} Z + \, \bar{\zeta} G_2,
\end{align}
where 
\begin{align}
\label{eq: hat check coeffs}
Z= \begin{pmatrix}
        Z_1 \\
        Z_2
    \end{pmatrix}, \quad 
    \bar{M}= (\bar \Sigma^*_{12})^\top (\bar\Sigma^*_{11})^{-1}, \quad 
    \bar{\zeta}^2 = \bar \Sigma^*_{22}  \, -  \, (\bar \Sigma^*_{12})^\top (\bar\Sigma^*_{11})^{-1} \bar \Sigma^*_{12}, \quad (G_1, G_2)^\top \sim N\left( 0, \begin{pmatrix}
        1 & q \\
        q &  1
    \end{pmatrix} \right),
\end{align}
and $(G_1, G_2)$ are independent of $Z$.
Similar to the proof strategy of \citet[Lemma~6.9]{donoho2016high}, we prove the following three claims:
\begin{enumerate}
    \item[(i)]   $|H(q)| \le H(1)=1$ for $q \in [0,1]$; 
    \item[(ii)] $H(q)$ is  strictly increasing for $q \in [0,1]$;
    \item[(iii)] $H(q)$ is convex for $q \in [0,1]$.
\end{enumerate}

To prove $(\rm i)$, we observe that at $q=1$, $H(1) =  (T^*_{B,11})^{-1} E[g(Z_1^*,   Z_1 + \sigma(Z_2) \bar\varepsilon; \, b_*)^2]$, which equals $1$, from the definition of $T^*_{B,11}$ in \eqref{eq: state evolution FP MT}. The inequality $|H(q)| \le H(1)$ for $q \in [0,1]$ follows from Cauchy-Schwarz. For $(\rm ii)$ and $(\rm iii)$, we use the same argument as  \citet[Lemma~C.1]{bayati2012lasso} and \citet[Appendix~3.2]{donoho2016high}.  
Using \eqref{eq: hat check exp}-\eqref{eq: hat check coeffs},  we define
\begin{align*}
    &h_{(Z, \bar\varepsilon)} (z) = g\left(\bar M Z + \bar \zeta z, \,   Z_1 + \sigma(Z_2) \bar\varepsilon; \, b_* \right),
\end{align*}
and write
\begin{align}
\label{eq: mcH def}
    H(q) = (T^*_{B,11})^{-1}E[\mathcal H (q; \, Z, \bar\varepsilon)], \qquad \mathcal H (q; \, Z, \bar\varepsilon) \coloneq E_q[ h_{(Z, \bar\varepsilon)} (G_1) h_{(Z, \bar\varepsilon)} (G_2) \mid Z, \bar \varepsilon],
\end{align}     where $E_q$ denotes the expectation with $(G_1, G_2)$ distributed as specified in \eqref{eq: hat check coeffs}. Letting $(X_t)_{t\ge 0}$ denote the standard Ornstein-Uhlenbeck process with covariance $E[X_0 X_t] = e^{-t}$, we have 
\[
\mathcal H (q; \, Z, \bar\varepsilon) = E[h_{(Z, \bar\varepsilon)} (X_0) h_{(Z, \bar\varepsilon)} (X_t)]\vert_{t = \log(1/q) },
\]
where the expectation is with respect to $X$. We now express $h_{(Z, \bar\varepsilon)}(\cdot) $ using the Hermite orthonormal basis $\{\varphi_l(\cdot) \}_{l\ge 0}$, and write $h_{(Z, \bar\varepsilon)} (X_0)= \sum_{l} c_l(Z, \bar\varepsilon) \varphi_l(X_0)$, where $c_l(Z, \bar\varepsilon) = E [h_{(Z, \bar\varepsilon)} (X_0) \varphi_l(X_0)]$, with the expectation taken over $X_0 \sim N(0,1)$. We recall that the Hermite polynomials $\varphi_l(\cdot)$ are eigenfunctions of the operator $P_t$ defined by $P_t f (x) = E[ f(X_t) \mid X_0 = x]$. Indeed, for $l \ge 0$, it holds that $P_t \varphi_l(x) = e^{-tl} \varphi_l(x)$ \citep{bakry2014analysis}. We then have the following spectral representation:
\begin{align*}
\mathcal H (q; \, Z, \bar\varepsilon) & = E [h_{(Z, \bar\varepsilon)} (X_0) \cdot P_t h_{(Z, \bar\varepsilon)} (X_0)]\vert_{t = \log (1/q)} \qquad &\\
  &=  E [\sum_{l = 0}^\infty c_l \varphi_l(X_0) P_t (\sum_{l' = 0}^\infty c_{l'}
  \varphi_{l'}(X_0))  ]\vert_{t = \log (1/q)}\\ 
 &=E [\sum_{l = 0}^\infty c_l \varphi_l(X_0) \sum_{l' = 0}^\infty c_{l'}  e^{-tl'}\varphi_{l'}(X_0)  ]\vert_{t = \log (1/q)} \\
  &= \sum_{l = 0}^\infty c_l^2 e^{-l \log(1/q) } = \sum_{l = 0}^\infty c_l^2 q^l.
\end{align*}
Claims $(\rm ii)$ and $(\rm iii)$ then follow by noting that $H(q) =  (T^*_{B,11})^{-1}\sum_{l = 0}^\infty E[c_l(Z, \bar\varepsilon)^2 ]q^l$, and  from the definition of $g$, it follows that $E[c_l(Z, \bar\varepsilon)^2 ] >0$ for at least one $l \ge 1$.

By convexity of $H(\cdot)$, we have 
\begin{equation}
\label{eq: Hq supporting line ineq}
    H(q) \ge H(1) - H'(1)(1-q) = 1 - H'(1)(1-q), \qquad q \in [0,1].
\end{equation} 
The proof is completed by showing that $0< H'(1) < 1$, which implies $H(q) > q$ for $q \in [0,1)$. Indeed, from \eqref{eq: mcH def} and using Price's theorem \citep{prices_theorem}, we obtain
\begin{align*}
    H'(q)  &  = (T^*_{B,11})^{-1} E \left[ E_q[ h'_{(Z, \bar\varepsilon)} (G_1) h'_{(Z, \bar\varepsilon)} (G_2) \mid Z, \bar \varepsilon] \right] \\
    & = (T^*_{B,11})^{-1} \bar{\zeta}^2 E \left[ E_q[ \partial_1 g\left(\bar M Z + \bar \zeta G_1, \,   Z_1 + \sigma(Z_2) \bar\varepsilon; \, b_* \right) \partial_1 g\left(\bar M Z + \bar \zeta G_2, \,   Z_1 + \sigma(Z_2) \bar\varepsilon; \, b_* \right) \mid Z, \bar \varepsilon]    \right],
\end{align*}
where the second equality follows from the definition of $h_{(Z, \bar\varepsilon)}$. Taking $q \to 1$, we obtain
\begin{equation*}
    \begin{split}
          H'(1) & =  (T^*_{B,11})^{-1} \bar{\zeta}^2E[\partial_1 g\left(\bar M Z + \bar \zeta G_2, \,   Z_1 + \sigma(Z_2) \bar\varepsilon; \, b_* \right)^2 ] \nonumber  \\
  & = (T^*_{B,11})^{-1} \bar{\zeta}^2 \frac{1}{1+b_*}E[-\partial_1 g\left(\bar M Z + \bar \zeta G_2, \,   Z_1 + \sigma(Z_2) \bar\varepsilon; \, b_* \right) ] \\
  & = (T^*_{B,11})^{-1} \bar{\zeta}^2 \frac{1}{(1+b_*)(\delta b_*)} \\  & = \frac{(\delta b_*) \bar{\zeta}^2}{(1+b_*) \tau_*^2}= \frac{b_*}{1+b_*}. 
    \end{split}
\end{equation*}
 The second equality above follows from \eqref{eq: g_1_derivatives}, noting that $(\partial_1 g)^2 = \frac{-1}{1+ b_*} \partial_1 g$; the third equality uses $\partial_1 g = -\partial_2 g$ and the definition of $b_*$ in \eqref{eq: fixed point state evolution}; the fourth equality uses  $\tau_*^2=T^*_{B,11} (\delta b_*)^2$; finally, the fifth holds since $\bar \zeta^2 = {\tau_*^2}/{\delta}$, which can be verified by direct calculation using the expressions for $\bar \Sigma_{12}^*$ and $\bar \Sigma_{22}^*$ in \eqref{eq: barSigma*_def}.
 
Since $b_* >0$, we have $H'(1) <1$, and thus by \eqref{eq: Hq supporting line ineq}, $H(q) > q$ for $q \in [0,1)$. Recalling from $(\rm i)$ that $H(q)$ is bounded,  monotone convergence implies $\lim_{k \to \infty}q_k =1$, which completes the proof. 
\end{proof}

\subsection{Proof of Lemma~\ref{lem: inliers near AMP iterate}}\label{ssec: appendix proof of sufficient inliers}

\begin{proof}
We first show $\lim_{\eta \to 0} \delta q_*(\eta) > 1 $. 
By the first equation of \eqref{eq: fixed point state evolution} and $g$ for Huber loss in \eqref{eq: g_1_derivatives}, for $b_* >0$, it holds that, 
\begin{eqnarray*}
    \lim_{\eta \to 0}  q_* (\eta) = P(|Z_1 + \sigma(Z_2)\bar\varepsilon - Z_1^* | < a(1 + b_*) ) = \delta^{-1}(1 + b_* ) /b_* =\delta^{-1}( 1 + b_*^{-1} )>  \delta^{-1}.
\end{eqnarray*}
Let $R_1^{k+1} = \bm Y - \Theta_1^{k+1}$ be as defined in Lemma~\ref{lem: vanishing increments} with components denoted by $R_{1,i}^{k+1}$. Then, using the Huber AMP update  in \eqref{eq: Huber_AMP} and recalling that $b_k = b_*$ for $k \ge 0$, we have
\begin{eqnarray}
 R_1^{k+1}
= \bm Y- \Theta_1^{k+1} 
= \bm Y-  \tilde{\bm X}\check\beta^{k+1}
+
b_* g(\Theta_1^k,\bm Y; b_*). \label{eq: prop proof R1 1}
\end{eqnarray}

Further, by \eqref{eq: effective score} 
\begin{align}\label{eq: prop proof proximal expansion}
\bar g(\Theta_1^k, \bm Y; b_*)
& =\Theta_1^k+b_* g(\Theta_1^k, \bm Y; b_*) & \nonumber
\\&=
\bm Y-R_1^k+b_* g(\Theta_1^k, \bm Y; b_*) 
= 
R_1^{k+1}-R_1^k + \tilde{\bm X}\check\beta^{k+1} \, ,
\end{align}
where the third equality uses \eqref{eq: prop proof R1 1}. 
Then, the ordinary residual can be expressed as
\begin{eqnarray}\label{eq: inlier residual expression 1}
\bm Y - \tilde{\bm X}\check\beta^{k+1} = (\bm Y - \bar g(\Theta_1^k, \bm Y; b_*) )+ (R_1^{k+1} - R_1^k) .   
\end{eqnarray}
For $i \in \{1, \ldots, n\}$, for the components satisfying $ |R_{1,i}^k| < (a- 2\eta) (1 + b_*)$, by \eqref{eq: effective score} and \eqref{eq: g_1 function Huber}, it follows that 
\begin{eqnarray}
| Y_i - \bar g(\Theta_{1i}^k, Y_i; b_*) | = |\frac{Y_i - \Theta_{1i}^{k}}{1 + b_*}| = 
   | \frac{R_{1,i}^k}{1 + b_*} |< (a-2\eta),
\end{eqnarray}
Thus, $| Y_i - \bar g(\Theta_{1i}^k, Y_i; b_*) | < (a- 2\eta)$ if $|R_{1,i}^k| <(a - 2\eta) (1 + b_*)$. 
Denoting $A_i \coloneq \{ |Y_i - \bar g(\Theta_{1i}^k, Y_i; b_*)| < (a- 2\eta) \}$ and $B_i \coloneq \{ |R_{1, i}^{k+1} - R_{1, i}^k | > \eta  \}$, it holds that 
\begin{align}\label{eq: inlier indicator inequality}
    I\{ | Y_i - \tilde X_i^\top \check\beta^{k+1} | <(a - \eta)\} &  \geq I\{ A_i \backslash B_i \}  \nonumber \\
    & \ge I\{ A_i\} - I\{B_i\} \ge I\{|R_{1,i}^k | < (a - 2\eta) (1 + b_*) \} -I\{B_i\} .
\end{align}
The second term, with probability 1, satisfies
\begin{eqnarray}\label{eq: inliear second term vanish}
\frac{1}{n}\sum_{i=1}^n I\left\{ |R_{1, i}^{k+1} - R_{1, i}^k| > \eta \right\}  \leq \frac{1}{\eta^2} \frac1n\|R_{1}^{k+1} - R_1^k\|_2^2 \stackrel{k\to\infty}{\to} 0 ,
\end{eqnarray}
where 
the inequality holds by Markov's inequality, and the convergence result is by \eqref{eq: increment vanish residual}.

By \eqref{eq: inlier indicator inequality} and \eqref{eq: inliear second term vanish}, we obtain
\begin{eqnarray}
\lefteqn{\lim\inf_{n,p\to \infty}q_k (\eta) = \frac{1}{n} \sum_{i=1}^n I\{ |Y_i - \tilde X_i^\top \check\beta^{k+1} | < (a - \eta)\}} &&\\
&\ge& \lim\inf_{n,p\to \infty}\frac{1}{n} \sum_{i=1}^n I\{ |R_{1,i}^k| <(a- 2\eta)(1 + b_*) \} - \lim\sup_{n,p\to \infty}\frac{1}{n}\sum_{i=1}^n I\{ |R_{1, i}^{k+1} - R_{1, i}^k| > \eta\}  \nonumber\\
&\ge &
\mathbb P(|Z_1 + \sigma(Z_2)\bar\varepsilon - Z_1^k| < (a - 2\eta ) (1 + b_*) ) -  
\lim_{n, p \to\infty}\frac{1}{\eta^2} \frac1n\|R_{1}^{k+1} - R_1^k\|_2^2 \nonumber \\
& \stackrel{k \to \infty}{\to}&  \mathbb P(|Z_1 + \sigma(Z_2)\bar\varepsilon - Z_1^*| < (a - 2\eta ) (1 + b_*) ) = q_*(\eta)  \nonumber,
\end{eqnarray}
where the convergence holds almost surely by Lemma~\ref{lem: wasserstein convergence of matrix AMP}.
\end{proof}

\subsection{Proof of Lemma \ref{lem: mse convergence Huber and AMP at stationarity}} \label{appendix: mse convergence Huber and AMP at stationarity}

\begin{proof}

Denote the objective in \eqref{eq: minimization whitened coordinate} and its gradient by
\begin{eqnarray}\label{eq: gradient definition}
  L_n(\beta)
= \sum_{i=1}^n \ell_a( \tilde X_i^\top \beta, Y_i), \qquad \qquad \nabla_\beta L_n(\beta)
=  \tilde{\bm X}^{\top}\partial_1\ell_a( \tilde{\bm X} \beta, \bm Y),
\end{eqnarray}
where $\partial_1\ell_a(\tilde{\bm X} \beta, \bm Y)$ is applied row-wise. Thus, $\check\beta$ is the minimizer of $L_n$.  By \eqref{eq: proximal estimating equations} and \eqref{eq: effective score},
\begin{eqnarray}\label{eq: prop proof proximal equivalence}
    \partial_1\ell_a(\bar g(u,y; b_*),y)
=
-g(u,y; b_*).
\end{eqnarray}
Recall that $R_1^{k+1} = \bm Y - \Theta_1^{k+1}$ is defined in Lemma~\ref{lem: vanishing increments} and $b_k = b_*$ for $k \ge 0$. Consider the second line of  the Huber AMP update \eqref{eq: Huber_AMP}. Recalling the stationarity of the state evolution and the fact that $\bar C_{11}^k = E[\partial_1 g(Z_1^*, \,  Z_1 + \sigma(Z_2) \bar\varepsilon; b_*)] = - (\delta b_* )^{-1}$,
we have
\begin{eqnarray}\label{eq: prop proof R1 2}
\check\beta^{k+1}-\check\beta^k
= \delta b_*
\tilde{\bm X}^{\top}g(\Theta_1^k, \bm Y; b_k).
\end{eqnarray}

Using the gradient expression in  \eqref{eq: gradient definition}, we have 
\begin{align}
& \frac{1}{\sqrt p}
\|\nabla_\beta L_n(\check\beta^{k+1})\|_2 = 
\frac{1}{\sqrt p}\|- \tilde{\bm X}^{\top}\partial_1\ell_a( \tilde{\bm X} \check\beta^{k+1},\bm Y) \|_2 \nonumber  \\
&= \frac{1}{\sqrt p}\|-  \tilde{\bm X}^{\top}\partial_1\ell_a( \tilde{\bm X} \check\beta^{k+1},\bm Y) +  \tilde{\bm X}^{\top}\partial_1\ell_a( \bar g(\Theta_1^k,\bm Y; b_k), \bm Y) -  \tilde{\bm X}^{\top}\partial_1\ell_a( \bar g(\Theta_1^k,\bm Y; b_k), \bm Y)\|_2
\nonumber
\\
&\leq \frac{1}{\sqrt p}\|\tilde{\bm X}\|_{\mathrm{op}} \|  \bar g(\Theta_1^k,\bm Y ; b_k)- \tilde{\bm X}\check\beta^{k+1}\|_2
 + 
 \frac{1}{\sqrt{p}}\|\tilde{\bm X}^{\top}\partial_1\ell_a( \bar g(\Theta_1^k,\bm Y; b_k), \bm Y)\|_2
 \nonumber
 \\
&\leq
\|\tilde{\bm X}\|_{\mathrm{op}}\sqrt{\frac{n}{p}}
\frac{\|R_1^{k+1}-R_1^k\|_2}{\sqrt n}
 + (\delta b_*)^{-1} \frac{\|\check\beta^{k+1}-\check\beta^k\|_2}{\sqrt{p}}
 ,
\label{eq: gradient control no psi}
\end{align}
where the first inequality uses the triangle inequality and the fact that $\partial_1 \ell_a$ is 1-Lipschitz in the first argument. The first term in the second inequality uses \eqref{eq: prop proof proximal expansion}, while the second term uses \eqref{eq: prop proof R1 2}.
By Lemma~\ref{lem: vanishing increments}, with probability one, we have 
\begin{eqnarray*}
 &&\lim_{k\to\infty}\lim_{n,p\to\infty}
\frac{1}{p}
\|\check\beta^{k+1}-\check\beta^k\|_2^2
=
0, \\
&&
\lim_{k\to\infty}\lim_{n,p\to\infty}
\frac{1}{n}
\|R_1^{k+1}-R_1^k\|_2^2
=
0.    
\end{eqnarray*}
Since  $\lim\sup_{n\to \infty}\|\tilde{\bm X}\|_{\rm op}=O(1)$ almost surely and $b_* >0$,
\eqref{eq: gradient control no psi} yields almost surely
\begin{equation}
\label{eq: lemma 3 gradient as to 0}
\lim_{k\to\infty}\limsup_{p\to\infty}
\frac{1}{\sqrt p}
\|\nabla_\beta L_n(\check\beta^{k+1})\|_2
=0.
\end{equation}

By Assumption~\ref{cond: restricted secant inequality} and \eqref{eq: lemma 3 gradient as to 0}, for sufficiently large $k$ and $n, p$ and $c_0 \in (0, \underline c)$ we have 
\begin{eqnarray}\label{eq: lemma 3 local gradient bound}
\kappa(\rho) \ge c_0, \qquad \qquad \|\nabla_\beta  L_n(\check\beta^{k+1})\|_2 < c_0 \rho \sqrt{p} .    
\end{eqnarray}

We first show by contradiction that any minimizer $\check\beta$ of $L_n$ satisfies 
$\|\check\beta - \check\beta^{k+1}\|_2 \leq \rho\sqrt{p}$. Assume $\|\check\beta - \check\beta^{k+1}\|_2 >\rho\sqrt{p}$, and define the unit vector $d=\frac{\check\beta - \check\beta^{k+1}}{\|\check\beta - \check\beta^{k+1}\|_2}$ and the boundary point $u = \check\beta^{k+1} + \rho\sqrt{p}d$.
By Assumption~\ref{cond: restricted secant inequality}, it holds that 
\begin{eqnarray}\label{eq: gradient inequality}
\rho \sqrt{p} d^\top ( \nabla_\beta L_n(u ) -\nabla_\beta L_n (\check\beta^{k+1})) & = \langle \nabla_\beta L_n (u ) -\nabla_\beta L_n(\check\beta^{k+1}), u - \check\beta^{k+1} \rangle \\ & \ge \kappa(\rho) \| u - \check\beta^{k+1}\|_2^2 \ge c_0 (\rho\sqrt{p})^2. \nonumber
\end{eqnarray}
Thus, \begin{eqnarray}\label{eq: gradient at boundary}
d^\top \nabla_\beta L_n(u ) \ge d^\top \nabla_\beta L_n(\check\beta^{k+1})  + c_0 \rho\sqrt{p} \ge c_0 \rho\sqrt{p} - \| \nabla_\beta L_n(\check\beta^{k+1}) \|_2 > 0,    
\end{eqnarray}
where the first inequality holds by \eqref{eq: gradient inequality}, the second by Cauchy--Schwarz inequality, and the last is by \eqref{eq: lemma 3 local gradient bound}.
Recall that $u$ is the boundary point along the ray from $\check\beta^{k+1}$ to $\check\beta$. Since $\|\check\beta - \check\beta^{k+1}\|_2 > \rho \sqrt{p}$, $\check\beta$ lies beyond $u$ on this ray. Since $L_n$ is convex, its directional derivative along $d$ is nondecreasing as we move away from $\check\beta^{k+1}$ along this ray. Hence, $d^\top\nabla L_n(\check\beta)
\geq d^\top\nabla L_n( u)$ for $\check\beta$ satisfying $\|\check\beta - \check\beta^{k+1}\|_2 > \rho \sqrt p$.
Since $\check\beta$ is a minimizer, $\nabla L_n(\check\beta) = 0$. Because $\|d\|_2=1$, by the Cauchy--Schwarz inequality, it holds that
\[
0 = d^\top\nabla L_n(\check\beta) \ge d^\top\nabla_\beta L_n( u)
> 0,
\]
a contradiction and the last inequality is by \eqref{eq: gradient at boundary}. Hence any minimizer  $\check \beta$ satisfies $\|\check\beta - \check\beta^{k+1}\|_2 \leq \rho\sqrt{p}$, i.e. it lies inside the ball of Assumption \ref{cond: restricted secant inequality}.

Taking $h= \check\beta - \check\beta^{k+1}$ in Assumption \ref{cond: restricted secant inequality}, we obtain that the following holds almost surely,  
\begin{eqnarray*}
\lim_{k\to\infty}\limsup_{p\to\infty}
\frac{
\|\check\beta^{k+1}-\check\beta\|_2
}{\sqrt p}
\leq
\frac{1}{c_0}
\lim_{k\to\infty}\limsup_{p\to\infty}
\frac{
\|\nabla_\beta L_n(\check\beta^{k+1})\|_2
}{\sqrt p}=0,
\end{eqnarray*}
where the equality follows from \eqref{eq: lemma 3 gradient as to 0}. 
\end{proof}

\section{Additional Figures}\label{appendix: figures}

\begin{figure}[t!]
    \centering
    \includegraphics[width=0.85\textwidth]{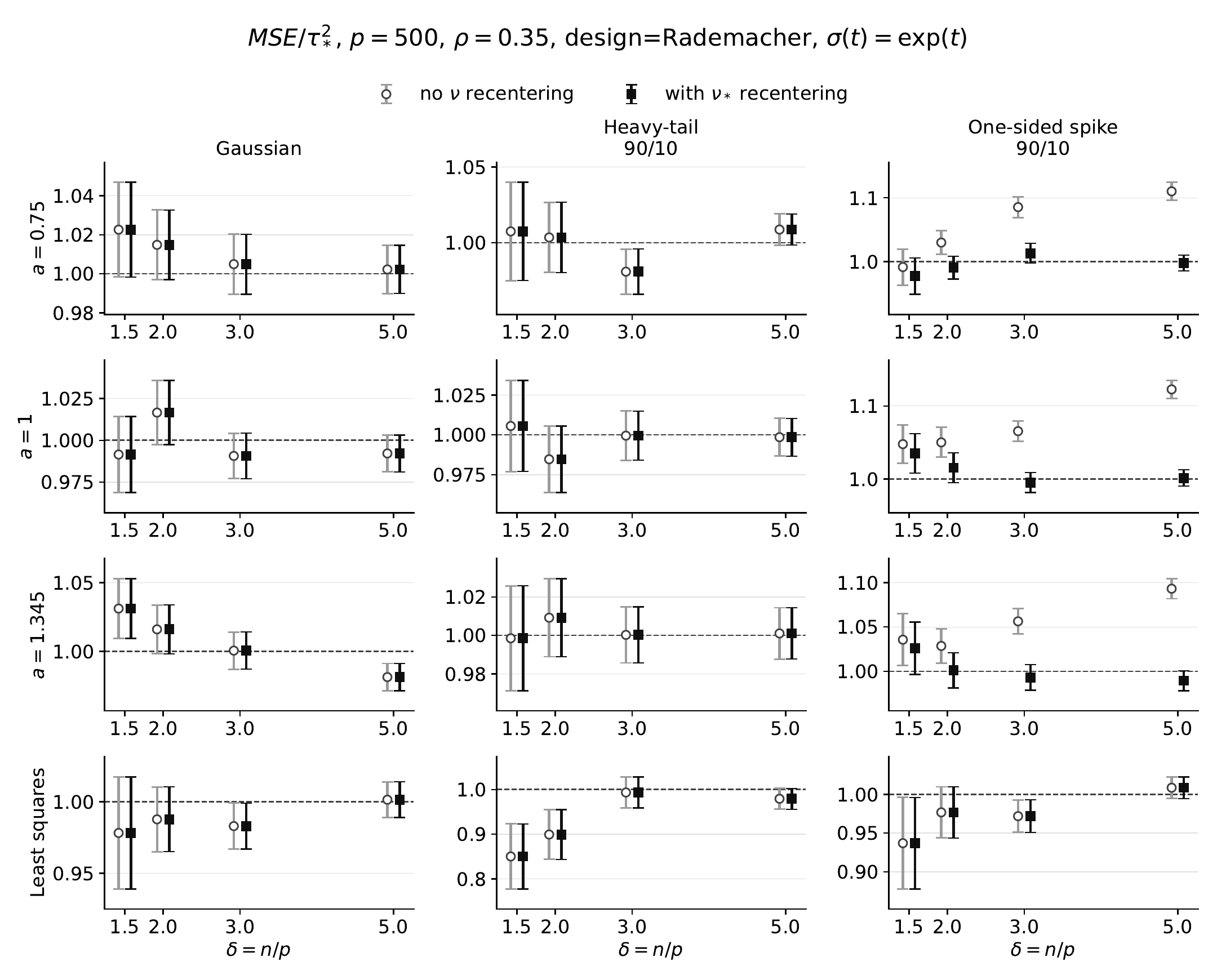}
    \caption{$Y_i = X_i^\top \beta_0 + \exp(X_i^\top \gamma_0)\varepsilon_i$, $X_i = \Omega_p^{1/2}\tilde X_i$,  $\mathbb P(\tilde X_{ij}=\frac{1}{\sqrt n})=\mathbb P(\tilde X_{ij}= \frac{-1}{\sqrt n})=\frac12$, and $(\Omega_p)_{jk}=0.35^{|j-k|}$. $p = 500$, $\delta = 1.5, 2, 3, 5$. The error $\varepsilon_i =  ( e_i - E(e_i))/\sqrt{\operatorname{var} (e_i)}$, where $e_i \sim N(0,1)$ (left column), $0.9N(0, 1) + 0.1 N(0, 16)$ (middle column), $0.9 N(0, 1) + 0.1 \delta_{8}$ (right column). The parameter vectors $\beta_0$ and $\gamma_0$ are generated once and held fixed across all replications. We set $\beta_{0,1} = 0$, $\gamma_{0,1} = 10$, $\beta_{0, 2} = 4$, and $\gamma_{0,2} = 2$. The other components of $\beta_0$ and $\gamma_0$ are randomly generated from a standard Gaussian distribution $N(0,1)$. The rows correspond to Huber threshold $a = 0.75, 1, 1.345, \infty$, respectively. The grey and black error bars (mean and 1.96 standard deviations over 200 replications) represent the ratios $R_{\rm uncentred}$ and $R_{\rm centred}$ in \eqref{eq: MSE ratios}. }
    \label{fig: recenter demo RademacherExp}
\end{figure}

\begin{figure}[t!]
    \centering
    \includegraphics[width=0.85\textwidth]{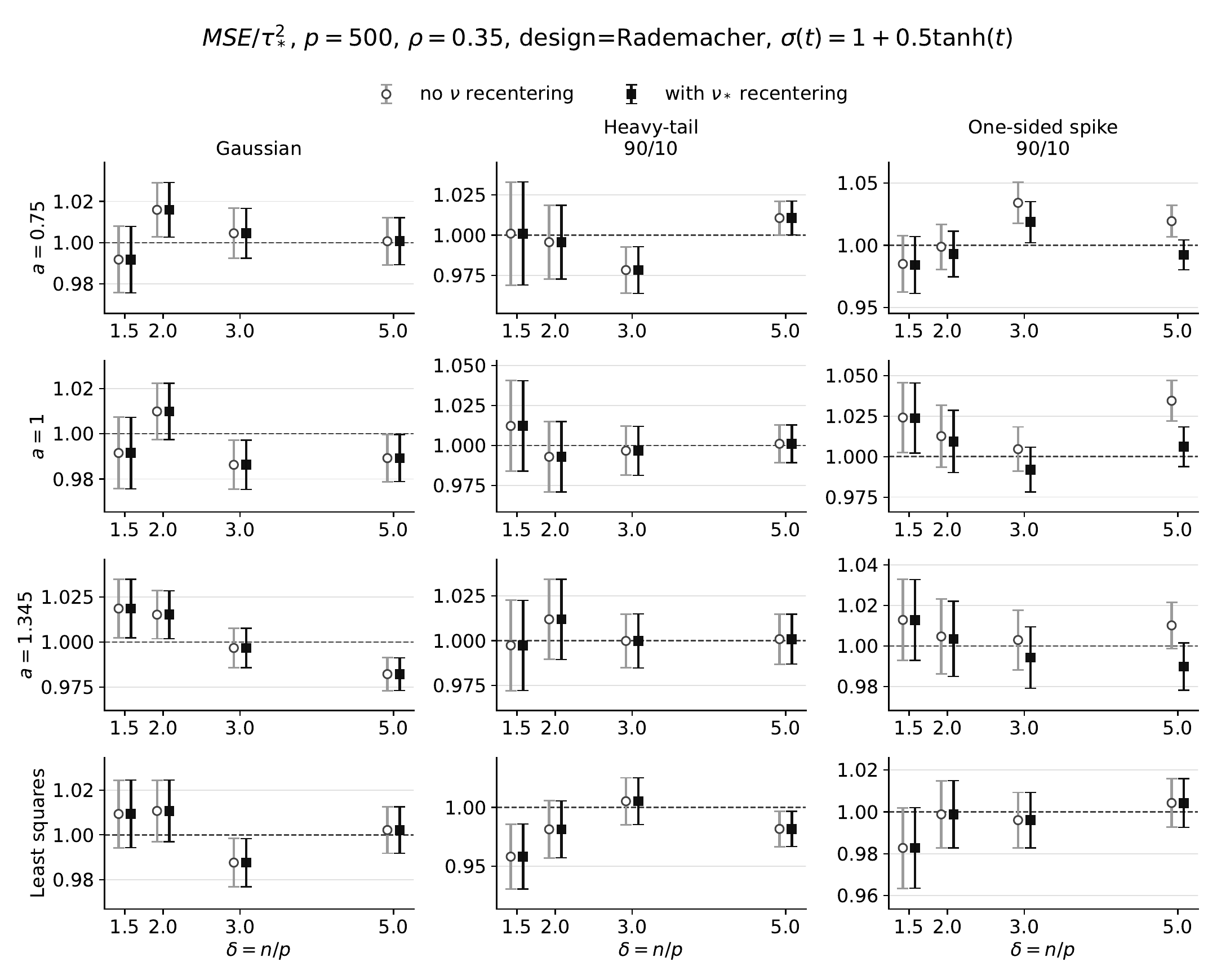}
    \caption{$Y_i = X_i^\top \beta_0 + (1 + 0.5\tanh(X_i^\top \gamma_0))\varepsilon_i$, $X_i = \Omega_p^{1/2}\tilde X_i$, $\mathbb P(\tilde X_{ij}=\frac{1}{\sqrt n})=\mathbb P(\tilde X_{ij}= \frac{-1}{\sqrt n})=\frac12$, and $(\Omega_p)_{jk}=0.35^{|j-k|}$. $p = 500$, $\delta = 1.5, 2, 3, 5$. The error $\varepsilon_i =  ( e_i - E(e_i))/\sqrt{\operatorname{var} (e_i)}$, where $e_i \sim N(0,1)$ (left column), $0.9N(0, 1) + 0.1 N(0, 16)$ (middle column), $0.9 N(0, 1) + 0.1 \delta_{8}$ (right column). The parameter vectors $\beta_0$ and $\gamma_0$ are generated once and held fixed across all replications. We set $\beta_{0,1} = 0$, $\gamma_{0,1} = 10$, $\beta_{0, 2} = 4$, and $\gamma_{0,2} = 2$. The other components of $\beta_0$ and $\gamma_0$ are randomly generated from a standard Gaussian distribution $N(0,1)$. The rows correspond to Huber threshold $a = 0.75, 1, 1.345, \infty$, respectively. The grey and black error bars (mean and 1.96 standard deviations over 200 replications) represent the ratios $R_{\rm uncentred}$ and $R_{\rm centred}$ in \eqref{eq: MSE ratios}.  }
    \label{fig: recenter demo RademacherTanh}
\end{figure}

\begin{figure}[t!]
    \centering
    \includegraphics[width=0.85\textwidth]{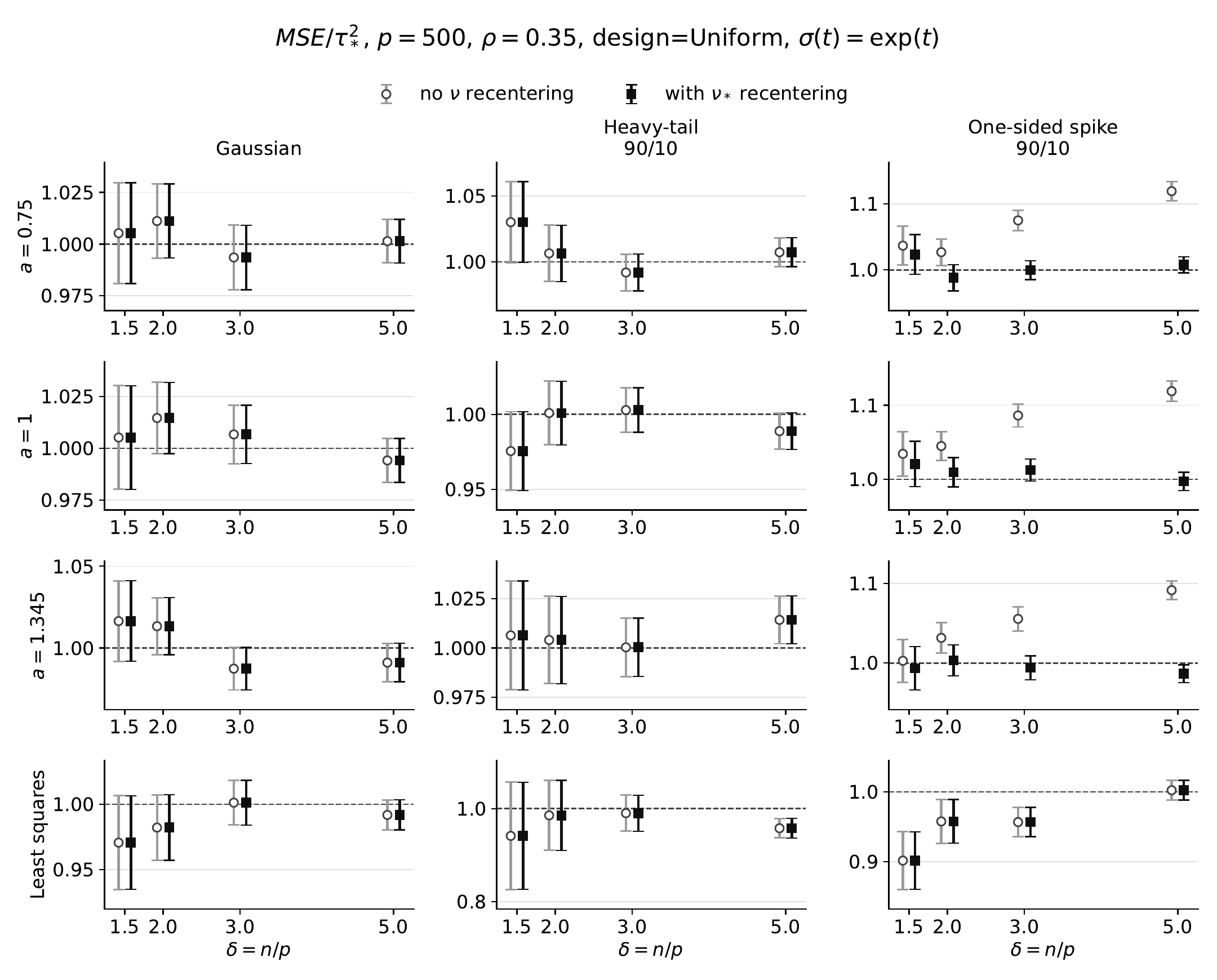}
    \caption{$Y_i = X_i^\top \beta_0 + \exp(X_i^\top \gamma_0)\varepsilon_i$, $X_i = \Omega_p^{1/2} \tilde X_i$, $\tilde X_{ij}\sim \operatorname{Unif}(-\sqrt{3/n},\sqrt{3/n})$, and $(\Omega_p)_{jk}=0.35^{|j-k|}$. $p = 500$, $\delta = 1.5, 2, 3, 5$. The error $\varepsilon_i =  ( e_i - E(e_i))/\sqrt{\operatorname{var} (e_i)}$, where $e_i \sim N(0,1)$ (left column), $0.9N(0, 1) + 0.1 N(0, 16)$ (middle column), $0.9 N(0, 1) + 0.1 \delta_{8}$ (right column). The parameter vectors $\beta_0$ and $\gamma_0$ are generated once and held fixed across all replications. We set $\beta_{0,1} = 0$, $\gamma_{0,1} = 10$, $\beta_{0, 2} = 4$, and $\gamma_{0,2} = 2$. The rest of the components of $\beta_0$ and $\gamma_0$ are randomly generated from a standard Gaussian distribution $N(0,1)$. The rows correspond to Huber threshold $a = 0.75, 1, 1.345, \infty$, respectively. The grey and black error bars (mean and 1.96 standard deviations over 200 replications) represent the ratios $R_{\rm uncentred}$ and $R_{\rm centred}$ in \eqref{eq: MSE ratios}. }
    \label{fig: recenter demo uniformExp}
\end{figure}

\begin{figure}[t!]
    \centering
    \includegraphics[width=0.85\textwidth]{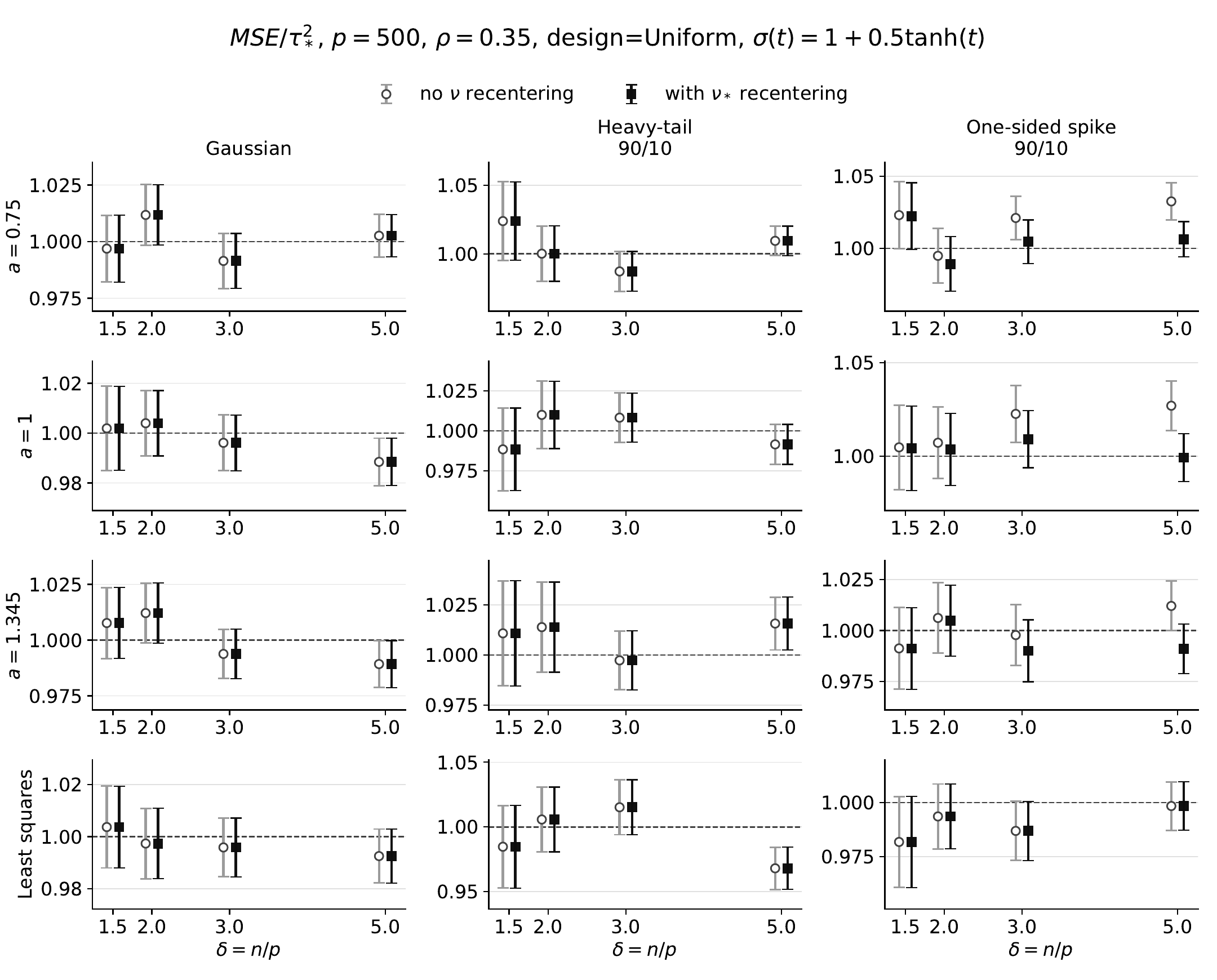}
    \caption{$Y_i = X_i^\top \beta_0 + (1 + 0.5\tanh(X_i^\top \gamma_0))\varepsilon_i$, $X_i = \Omega_p^{1/2}\tilde X_i$, $\tilde X_{ij}\sim \operatorname{Unif}(-\sqrt{3/n},\sqrt{3/n})$, and $(\Omega_p)_{jk}=0.35^{|j-k|}$. $p = 500$, $\delta = 1.5, 2, 3, 5$. The error $\varepsilon_i =  ( e_i - E(e_i))/\sqrt{\operatorname{var} (e_i)}$, where $e_i \sim N(0,1)$ (left column), $0.9N(0, 1) + 0.1 N(0, 16)$ (middle column), $0.9 N(0, 1) + 0.1 \delta_{8}$ (right column). The parameter vectors $\beta_0$ and $\gamma_0$ are generated once and held fixed across all replications. We set $\beta_{0,1} = 0$, $\gamma_{0,1} = 10$, $\beta_{0, 2} = 4$, and $\gamma_{0,2} = 2$. The rest of the components of $\beta_0$ and $\gamma_0$ are randomly generated from a standard Gaussian distribution $N(0,1)$. The rows correspond to Huber threshold $a = 0.75, 1, 1.345, \infty$, respectively. The grey and black error bars (mean and 1.96 standard deviations over 200 replications) represent the ratios $R_{\rm uncentred}$ and $R_{\rm centred}$ in \eqref{eq: MSE ratios}.}
    \label{fig: recenter demo uniformTanh}
\end{figure}

In this section, we present results from simulations where the whitened covariates are drawn from either a Rademacher (Figures \ref{fig: recenter demo RademacherExp}, \ref{fig: recenter demo RademacherTanh}) or a uniform distribution (Figures \ref{fig: recenter demo uniformExp}, \ref{fig: recenter demo uniformTanh}). The remaining parameters of the simulations are identical to those used for the Gaussian design simulations in Figures \ref{fig: gaussian exp} and \ref{fig: gaussian tanh}.  The results are qualitatively identical to those for Gaussian designs, consistent with the fact that Theorem \ref{thm: dist convergence of huber mean} holds for any whitened covariate distribution satisfying Assumption \ref{cond:white_noise_design}.

\end{document}